\documentclass[12pt,reqno]{amsart}
\usepackage[margin=1in]{geometry}
\usepackage{amsmath,amssymb,amsthm,mathtools}
\usepackage{graphicx}
\usepackage{booktabs}
\usepackage{float}
\usepackage{enumitem}
\usepackage{tikz}
\usepackage{xcolor}
\usepackage{hyperref}
\usepackage{microtype}

\newcommand{\red}[1]{\textcolor{red}{#1}}
\hypersetup{colorlinks=true,linkcolor=blue,citecolor=blue,urlcolor=blue}
\allowdisplaybreaks

\newtheorem{theorem}{Theorem}[section]
\newtheorem{proposition}[theorem]{Proposition}
\newtheorem{lemma}[theorem]{Lemma}

\theoremstyle{remark}
\newtheorem{remark}[theorem]{Remark}

\newcommand{\E}{\mathbb E}
\newcommand{\Pp}{\mathbb P}
\newcommand{\Var}{\operatorname{Var}}
\newcommand{\Cov}{\operatorname{Cov}}
\newcommand{\1}{\mathbf 1}
\newcommand{\N}{\mathbb N}

\newcommand{\calF}{\mathcal F}
\newcommand{\calP}{\mathcal P}
\newcommand{\LCM}{\operatorname{lcm}}
\newcommand{\dd}{\,\mathrm d}

\newcommand{\toP}{\xrightarrow{\Pp}}
\newcommand{\tod}{\xrightarrow{d}}

\makeatletter
\@namedef{subjclassname@2020}{\textup{2020} Mathematics Subject Classification}

\title[Phase Transitions for Random Least Common Multiples]
{Phase Transitions for Random Least Common Multiples with Memory}

\author{Shaochen Wang}
\address[S. Wang]{School of Mathematics, South China University of Technology, Guangzhou 510640, China}
\email{mascwang@scut.edu.cn}

\author{Guangyu Yang}
\address[G. Yang]{School of Mathematics and Statistics, Zhengzhou University, Zhengzhou 450001, China}
\email{guangyu@edu.zzu.cn}

\author{Wang Zhou}
\address[W. Zhou]{Department of Statistics and Data Science, National University of Singapore, Singapore 119077, Singapore}
\email{stazw@nus.edu.sg}

\date{\today}



\begin{document}

\begin{abstract}
We study the least common multiple of a random subset of
$\{1,\ldots,n\}$ when the selection follows a correlated Bernoulli process with memory.
The copying rule preserves the marginal retention
probabilities and creates a genealogical forest. At fixed
density, the logarithmic least common multiple exhibits a fluctuation
phase transition at the memory parameter $a=1/2$. The boundary limits
depend on how genealogical families are distributed among the divisibility
classes carrying the leading arithmetic weight. Our proofs combine the
von Mangoldt representation with a reduction to large primes. At fixed
density, an exact expansion separates an innovation martingale from a
memory term. An asymptotic projection onto the retained count identifies
the contribution of memory. Marked point-process convergence and estimates
for rare genealogical families give the boundary limits.
\end{abstract}

\maketitle

\setcounter{tocdepth}{1}
\tableofcontents


\section{Introduction}\label{sec:introduction}


Let $A_n$ be a random subset of $[n]:=\{1,\ldots,n\}$, and let
$L_n:=\LCM(A_n)$, with $\LCM(\varnothing)=1$. The least common multiple
(LCM) depends on both the number of selected integers and their divisibility
properties. For every finite set $A$ of positive integers,
\begin{equation}\label{eq:intro-vm}
   \log \LCM(A)=\sum_{m\ge1}\Lambda(m)
   \1_{\{A\cap m\N\ne\varnothing\}}.
\end{equation}
Here $m\N$ denotes the positive multiples of $m$, and the von Mangoldt
function satisfies $\Lambda(p^r)=\log p$ for primes $p$ and integers $r\ge1$,
and vanishes otherwise. Thus $\log L_n$ is a weighted sum of indicators that
a prime-power divisibility class contains a selected integer.

For independently retained integers, Cilleruelo, Ru\'e, \v{S}arka and
Zumalac\'arregui established the first asymptotic results for random LCMs
\cite{CRSZ}. Alsmeyer, Kabluchko and Marynych subsequently proved a strong
law, endpoint and functional central limit theorems, and two boundary
Poisson limits \cite{AKM}. At a fixed retention probability
$\theta\in(0,1)$, their endpoint theorem states that
\begin{equation}\label{eq:logLnCLT}
   \frac{\log L_n-\E\log L_n}{\sqrt{n\log n}}
   \tod N\bigl(0,g(1-\theta)\bigr),
   \qquad
   g(z)=\sum_{k\ge1}\frac{z^k(1-z^k)}{k(k+1)}.
\end{equation}
Functional large and moderate deviation principles and a functional law
of the iterated logarithm were obtained by Wang, Yang and Zhou in \cite{WYZ}. Related LCM models
based on uniform samples, growing tuples, and other random sets have been
studied in \cite{BostanMarynychRaschel,BuraczewskiIksanovMarynych,Sanna}.

A key step in the independent model is the reduction to primes larger
than $\sqrt n$. Two such primes cannot divide the same integer in $[n]$.
Their sets of multiples are disjoint, so independent retention makes the
corresponding indicators independent. This raises a question about the role
of the selection rule: {\it do the same limits hold when the marginal
retention probabilities are preserved but the decisions become dependent?} We study this question through a correlated Bernoulli process with
memory. For each $n$, let $X_{n,j}$ be the retention indicator at step $j$.
The first indicator has the Bernoulli$(\theta_n)$ distribution. At time
$j+1$, the new indicator copies a uniformly chosen
predecessor with probability $a\in[0,1)$. Otherwise, it is sampled afresh
from the same Bernoulli law. Every indicator has marginal retention
probability $\theta_n$, and $a=0$ gives independent retention. Thus $\theta_n$ controls the marginal density, while $a$ controls the dependence. When
$a>0$, copying can connect decisions at indices separated by a distance of order $n$. Even
for large primes with disjoint sets of multiples, the divisibility
indicators may then be dependent.

The same selection rule provides a model for the common period of randomly
activated components. Give component $j$ the phase $t\bmod j$ at integer
time $t$, and activate it when $X_{n,j}=1$. The first positive time at which
all activated components return to their initial state is $L_n$. In periodic
task systems, this common period is called the hyperperiod. Its length is
relevant to cyclic schedules and the generation of synthetic workloads
\cite{RipollBallesterRipoll}. The same divisibility condition makes $L_n$
the exponent of $\prod_{j\in A_n}\mathbb Z/j\mathbb Z$ and the order of a
permutation whose disjoint cycles have lengths $j\in A_n$. In this interpretation,
each component keeps its prescribed period and may
inherit the activation decision of an earlier component. The model can
therefore be used to generate periodic configurations with a fixed marginal
inclusion probability and varying dependence. The limit laws for $\log L_n$
describe how that dependence affects the common period on a logarithmic
scale.

\begin{figure}[H]
\centering
\begin{tikzpicture}[
vertex/.style={circle,draw,minimum size=7mm,inner sep=0pt},
positive/.style={vertex,fill=blue!12},
negative/.style={vertex,fill=gray!25}
]
\node[positive] (r1) at (0,1.2) {$1$};
\node[positive] (v5) at (1.4,1.8) {$5$};
\node[positive] (v6) at (1.4,0.6) {$6$};
\node[negative] (r2) at (3.0,1.2) {$2$};
\node[negative] (r3) at (4.6,1.2) {$3$};
\node[negative] (v7) at (6.0,1.2) {$7$};
\node[positive] (r4) at (7.6,1.2) {$4$};
\node[negative] (r8) at (12.2,1.2) {$8$};
\node[positive] (v9) at (10.6,1.8) {$9$};
\node[positive] (v10) at (10.6,0.6) {$10$};
\draw(v5)--(r1);
\draw (v6)--(r1);
\draw(v7)--(r3);
\draw (v9)--(r4);
\draw (v10)--(r4);

\end{tikzpicture}
\caption{\scriptsize An example of a genealogical forest with $10$ vertices. Each edge
joins a copied decision to its chosen predecessor. Vertices in the same
component inherit the label of their root, whereas fresh decisions start new
components. Blue vertices belong to components whose root label is $1$, and
gray vertices belong to components whose root label is $0$.}
\label{fig:genealogical-forest}
\end{figure}
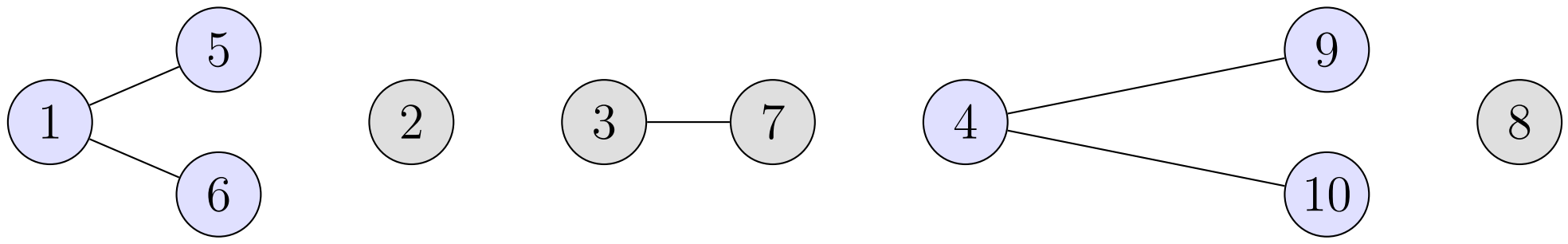

The underlying correlated Bernoulli process was introduced as a generalized
binomial construction in \cite{DreznerFarnum}. Its partial sums have been
studied by martingale methods in \cite{Heyde,JamesJamesQi,GonzalezLambert}.
They also form a positively step-reinforced random walk with Bernoulli
innovations, placing the model within the theory of elephant random walks
and related processes with memory; see
\cite{SchutzTrimper,HuZhang,ZhangStochastic,RoyTakeiTanemura,
RoyTakeiTanemuraVarying,DhillonKataria,NakanoTakei,
GuerinLaulinRaschel,HuDong,AguechQin,QinRecurrence,
PeresQinTransition,PeresQinMixing,
AguechBenHarizElMachkouriFaouzi,GutStadtmuller,FangTraining}.
For the LCM, the additional problem is to relate this temporal dependence
to the divisibility indicators in \eqref{eq:intro-vm}. The copying rule provides a genealogical representation of the selected
set. Join each copied vertex to its chosen predecessor, and let every fresh
decision start a new component. Assign independent Bernoulli$(\theta_n)$
labels to the roots and pass each root label to all its descendants. The
selected set is then the union of the components with label one in a
percolated random recursive tree. Common ancestry now determines which decisions are linked. To study
$L_n$, one needs to know how these components are distributed among the
prime-power divisibility classes.
Figure~\ref{fig:genealogical-forest} illustrates the construction.

Our results cover the fixed-density regime as well as the two boundary regimes.
At fixed density, a transition occurs at $a=1/2$, when fluctuations
of the retained count reach the arithmetic scale. At the boundaries, the
limiting laws depend on how rare families are distributed among the relevant
divisibility classes. Section~\ref{sec:model-main} gives the precise statements. Our proofs
combine arithmetic estimates with the genealogical representation.
At fixed density we pass to primes larger than $\sqrt{n}$, whose disjoint multiple sets allow an exact expansion that splits the fluctuations into a martingale of fresh decisions and a memory-carrying drift. An asymptotic projection onto
the retained count determines the contribution of memory. At the boundaries we show that, after the early roots are truncated, the rare positive and deletion families converge to a marked Poisson process.

The paper is organized as follows. Section~\ref{sec:model-main} introduces the
model and states the main results. Section~\ref{sec:fixed-proof} proves the fixed-density theorem,
and Section~\ref{sec:variable-proof} treats the two boundary regimes.
Section~\ref{sec-final} contains concluding remarks and further directions.
Appendix~\ref{app:auxiliary-estimates} collects the analytic and probabilistic
results used in the proofs.


\section{Model, notation, and main results}\label{sec:model-main}


Throughout the paper, $\calP$ is the set of all primes. We write
$\N=\{1,2,\ldots\}$ and use natural logarithms. Empty sums and products
are zero and one, respectively. Constants denoted by $C$ may change from
line to line; they are independent of $n$, and any additional dependence
needed in a uniform estimate is indicated. Unless stated otherwise,
$p$ denotes a prime, while $i,j,k,m,n$ denote positive integers. We write
$\toP$ and $\tod$ for convergence in probability and in distribution,
respectively. For $r>0$ and a random variable $Y$ with $\E|Y|^r<\infty$, set
$\|Y\|_r=(\E|Y|^r)^{1/r}$. We write $\Pi(\mu)$ for a Poisson random
variable with mean $\mu$. The notation $o_{L^r}(b_n)$ means convergence to
zero in $L^r$ after division by $b_n$; the corresponding probabilistic
notation is
$o_P(b_n)$. We use
\[
   \pi(x):=\#\{p\in\calP:p\le x\},
   \qquad
   \vartheta(x):=\sum_{p\le x}\log p
\]
for the prime counting function and the first Chebyshev function. The von
Mangoldt function is
\[
   \Lambda(m):=
   \begin{cases}
      \log p,&m=p^r\text{ for some }p\in\calP\text{ and }r\ge1,\\
      0,&\text{otherwise},
   \end{cases}
\]
and
\[
   \psi(x):=\sum_{m\le x}\Lambda(m)
\]
is the second Chebyshev function; see \cite[Chapters~2 and~4]{Apostol} for
these standard arithmetic functions. Prime number estimates used below are collected in
Appendix~\ref{app:auxiliary-estimates}.

For each $n\ge1$, let $\theta_n\in(0,1)$ and define the correlated Bernoulli sequence
\[
   X_{n,1},X_{n,2},\ldots,X_{n,n}\in\{0,1\}
\]
as follows. First, $X_{n,1}$ is Bernoulli with parameter $\theta_n$. Given
$\calF_{n,j}:=\sigma(X_{n,1},\ldots,X_{n,j})$, set
\begin{equation}\label{eq:universal-model}
   \Pp(X_{n,j+1}=1\mid\calF_{n,j})
   =\theta_n+a\bigg(\frac1j\sum_{i=1}^jX_{n,i}-\theta_n\bigg),
   \qquad 1\le j<n,
\end{equation}
where $0\le a<1$. The memory parameter $a$ is fixed throughout every
limit. By the recursion, we have $\E X_{n,j}=\theta_n$ for every $j\le n$.
Equivalently, with probability $a$, the new variable copies one of the
preceding variables chosen uniformly. With probability $1-a$, it
is sampled independently from the Bernoulli$(\theta_n)$ law. Define
\[
   A_n:=\{j\le n:X_{n,j}=1\},
   \qquad
   L_n:=\LCM(A_n),
\]
with $\LCM(\varnothing)=1$, and denote
\[
   T_n:=\log L_n,
   \qquad
   M_n:=|A_n|.
\]
Thus $M_n$ is the number of retained integers and $T_n$ is the logarithmic
least common multiple.

We call the component construction described in the introduction the
\emph{genealogical forest}. Conditional on this forest, the root labels are
independent Bernoulli$(\theta_n)$ variables, and every vertex has the label of
its root.


\subsection{Fixed density}\label{sec:fixed-result}


Assume that $\theta_n\equiv\theta\in(0,1)$ and write $q=1-\theta$.
For almost sure and $L^2$ limits, take $X_{n,j}=X_j$, $j\le n$, from one
infinite sequence satisfying \eqref{eq:universal-model}.
Set
\[
D_n:=M_n-n\theta.
\]
For $1/2<a<1$, let $W_{a,\theta}$ denote the nondegenerate $L^2$
limit of $D_n/n^a$; it has mean zero, and its existence is stated in
Proposition~\ref{prop:count-transition}. For $1/2\le a<1$, define
\begin{equation}\label{eq:Kdef}
   K_{a,\theta}:=
   \sum_{k\ge1}q^{k-1}\bigl(k^{-a}-(k+1)^{-a}\bigr)
   \sum_{r=1}^k r^{a-1}.
\end{equation}
This series is strictly positive and converges absolutely because its
$k$th summand is $O(q^k/k)$.

\begin{theorem}\label{thm:fixed-main}
Let $\theta_n\equiv\theta\in(0,1)$ and $0\le a<1$. Then
\begin{equation}\label{eq:main-lln}
   \frac{T_n}{n}\toP
   c(\theta):=\frac{\theta\log(1/\theta)}{1-\theta}.
\end{equation}
For the centered fluctuations, the following limits hold.
\begin{enumerate}[label=(\roman*)]
\item If $0\le a<1/2$, then
\begin{equation}\label{eq:main-weak}
   \frac{T_n-\E T_n}{\sqrt{n\log n}}
   \tod N\bigl(0,g(q)\bigr),
   \qquad
   g(q)=\sum_{k\ge1}\frac{q^k(1-q^k)}{k(k+1)}.
\end{equation}
\item If $a=1/2$, then
\begin{equation}\label{eq:main-critical}
   \frac{T_n-\E T_n}{\sqrt{n\log n}}
   \tod N\bigl(0,g(q)+\theta qK_{1/2,\theta}^2\bigr).
\end{equation}
\item If $1/2<a<1$, then
\begin{equation}\label{eq:main-strong}
   \frac{T_n-\E T_n}{n^a}
   \tod K_{a,\theta}W_{a,\theta}.
\end{equation}
\end{enumerate}
\end{theorem}

When $a=0$, Theorem \ref{thm:fixed-main} reduces to the endpoint theorem of
Alsmeyer, Kabluchko and Marynych \cite{AKM}. Throughout the weak-memory range
$a<1/2$, both the normalization and the limiting variance are the same as
under independent retention. At $a=1/2$, dependence first appears at the
leading arithmetic scale through an additional,
asymptotically independent Gaussian count term. For $a>1/2$, the scale
changes from $\sqrt{n\log n}$ to $n^a$, and the non-Gaussian count process
limit determines the leading fluctuation. The law of large numbers remains
unchanged in all three regimes.


\subsection{Sparse retention}\label{subsec:sparse-retention}


We next consider the boundary regime
\begin{equation}\label{eq:sparse-assumption}
   \theta_n\sim\frac{\lambda}{n},
   \qquad \lambda>0.
\end{equation}
In this regime, $\E M_n=n\theta_n\to\lambda$, so the retained count remains
of order one. For $a>0$, dependence remains visible because one positive root
may generate several retained descendants. Indeed, a vertex
$j\ge2$ starts a new component with probability $1-a$, and this component
receives a positive label with probability $\theta_n$.
Thus positive roots occur with probability $(1-a)\theta_n\sim(1-a)\lambda/n$.
After rescaling their locations by $n$, these rare roots therefore converge
to a Poisson point process $\mathcal R$ on $(0,1]$ with intensity
$(1-a)\lambda\,\dd u$. A root born near time $un$ may subsequently be copied.
Conditional on $\mathcal R$, assign independent marks $Z_u$ to its
points. For $a>0$, the mark at $u$ has the geometric distribution
\begin{equation}\label{eq:Yule-mark}
   \Pp(Z_u=m)
   =u^a(1-u^a)^{m-1},
   \qquad m\ge1.
\end{equation}
For $a=0$, we set $Z_u\equiv1$.
The limiting count is the total mass of the marked Poisson process,
\begin{equation}\label{eq:sparse-limit-var}
   N_{a,\lambda}:=\sum_{u\in\mathcal R}Z_u.
\end{equation}
$N_{a,\lambda}$ also has a compound Poisson representation. For $a>0$, let
$U$ be uniformly distributed on $(0,1)$ and, conditional on $U=u$, let $Z$
have the distribution in \eqref{eq:Yule-mark}. By mixing over $U$, we obtain
\begin{equation}\label{eq:Yule-Simon-mark}
   \Pp(Z=m)
   =\int_0^1u^a(1-u^a)^{m-1}\,\dd u
   =\frac1a B\Big(m,1+\frac1a\Big),
   \qquad m\ge1,
\end{equation}
where $B$ is the beta function. Thus $Z$ has the Yule--Simon distribution
with parameter $1/a$ \cite{Simon,BaurBertoinYule}. If $J\sim\operatorname{Poisson}((1-a)\lambda)$ and $Z_1,Z_2,\ldots$ are
independent copies of $Z$, independent of $J$, then
\begin{equation}\label{eq:sparse-compound-poisson}
   N_{a,\lambda}
   \ \stackrel{d}{=}\
   \sum_{k=1}^{J}Z_k.
\end{equation}
Equivalently, its probability generating function is
\begin{equation}\label{eq:sparse-pgf}
   \E[z^{N_{a,\lambda}}]
   =\exp\left\{(1-a)\lambda\int_0^1
   \left(
      \frac{u^az}{1-(1-u^a)z}-1
   \right)\dd u\right\},
   \qquad 0\le z\le1.
\end{equation}

\begin{theorem}\label{thm:sparse-main}
Let $0\le a<1$ and assume \eqref{eq:sparse-assumption}. Then
\begin{equation}\label{eq:sparse-count}
   M_n\tod N_{a,\lambda},
\end{equation}
and
\begin{equation}\label{eq:sparse-lcm}
   \frac{T_n}{\log n}\tod N_{a,\lambda}.
\end{equation}
\end{theorem}

The conditional mark moments are
\begin{equation}\label{eq:Yule-first-two-moments}
   \E[Z_u]=u^{-a},
   \qquad
   \E[Z_u^2]=2u^{-2a}-u^{-a}.
\end{equation}
Hence $\E N_{a,\lambda}=\lambda$, and if $a<1/2$, then
\begin{equation}\label{eq:sparse-variance}
   \Var(N_{a,\lambda})
   =\frac{\lambda}{1-2a},
\end{equation}
whereas $N_{a,\lambda}$ has an infinite second moment when $a\ge1/2$.
For $a=0$, Theorem~\ref{thm:sparse-main} reduces to the independent Poisson
limit in \cite[Theorem~1.7]{AKM}. In this case every mark equals one, so
$N_{0,\lambda}$ is Poisson with mean $\lambda$.


\subsection{Near-complete retention}\label{subsec:near-complete-retention}


When almost all integers are retained, we study the logarithmic defect
\[
   \Delta_n:=\psi(n)-T_n.
\]
Thus $\Delta_n$ is the logarithmic loss relative to the least common
multiple of all integers in $[n]$. Under \eqref{eq:nearcomplete-assumption},
the leading defect comes from missing primes in $(n/2,n]$ and has the same
Poisson limit as under independent retention.

\begin{theorem}\label{thm:near-complete}
Let $0\le a<1$ and assume
\begin{equation}\label{eq:nearcomplete-assumption}
   1-\theta_n\sim\frac{\lambda\log n}{n},
\end{equation}
with fixed $\lambda>0$. Then
\begin{equation}\label{eq:nearcomplete-limit}
   \frac{\Delta_n}{\log n}\tod\Pi(\lambda/2).
\end{equation}
\end{theorem}

For $a=0$, Theorem~\ref{thm:near-complete} reduces to the independent
boundary law in \cite[Theorem~1.8]{AKM}.


\section{Proof of the fixed-density theorem}\label{sec:fixed-proof}


Throughout this section, $\theta_n\equiv\theta\in(0,1)$ and $q=1-\theta$.
We realize all rows as initial segments of one infinite correlated Bernoulli
process with parameters $a$ and $\theta$, and therefore omit the row index.
Thus $\calF_j=\sigma(X_1,\ldots,X_j)$.
Denote
\[
D_n=M_n-n\theta=\sum_{j=1}^n(X_j-\theta),
\]
and
\[
\gamma_1:=1,
\qquad
\gamma_n:=\prod_{j=1}^{n-1}\Big(1+\frac aj\Big)
=\frac{\Gamma(n+a)}{\Gamma(1+a)\Gamma(n)}.
\]
Here $\Gamma$ is Euler's gamma function, and
$\gamma_n\sim n^a/\Gamma(1+a)$. We first recall the fluctuation theory of
$D_n$ and then combine genealogical estimates with the arithmetic
decomposition of $T_n$.


\subsection{Count process and genealogy}


The following proposition summarizes the count fluctuations of the
correlated Bernoulli process when $\theta$ is fixed.

\begin{proposition}
\label{prop:count-transition}
Let $0\le a<1$. Then $\mathcal M_n:=D_n/\gamma_n$ is a square-integrable
martingale. As $n\to\infty$, the following three regimes hold:
\begin{align}
0\le a<\frac12:
&\qquad
   \frac{D_n}{\sqrt n}\tod N\bigg(0,\frac{\theta q}{1-2a}\bigg),\label{eq:D-weak}\\
 a=\frac12:
&\qquad
   \frac{D_n}{\sqrt{n\log n}}\tod N(0,\theta q),\label{eq:D-critical}\\
\frac12<a<1:
&\qquad
   \frac{D_n}{n^a}\longrightarrow W_{a,\theta}
   \quad\text{almost surely and in }L^2.\label{eq:D-strong}
\end{align}
The last limit is nondegenerate, non-Gaussian, and has mean zero.
Moreover, $D_n/n\to0$ almost surely.
\end{proposition}

\begin{proof}
The convergence assertions are the constant parameter cases of \cite[Theorems~2.1 and~2.2]{JamesJamesQi}; see also \cite{Heyde}. In particular, the $L^2$ convergence follows from \cite[equation~(3.4)]{JamesJamesQi}, and the non-Gaussian assertion is the constant parameter result of \cite{Heyde}, also stated explicitly in \cite[p.~2339]{JamesJamesQi}.
\end{proof}

For later use, let $\calF_0$ be the trivial $\sigma$-field, set
$p_1=\theta$, $\xi_1=X_1-\theta$, and, for $j\ge2$, define
\[
   p_j:=\E(X_j\mid\calF_{j-1})
      =\theta+\frac{aD_{j-1}}{j-1},
   \qquad
   \xi_j:=X_j-p_j.
\]
Then $\E(\xi_j\mid\calF_{j-1})=0$, $|\xi_j|\le1$, and
\begin{equation}\label{eq:count-recursion-detailed}
   D_j=\left(1+\frac{a}{j-1}\right)D_{j-1}+\xi_j,
   \qquad j\ge2.
\end{equation}
By simple iteration, we have the martingale representation
\begin{equation}\label{eq:count-martingale-representation}
   \mathcal M_n:=\frac{D_n}{\gamma_n}
      =\sum_{j=1}^n\frac{\xi_j}{\gamma_j}.
\end{equation}

\begin{lemma}\label{lem:count-moments-localization}
For every $r\ge2$, there is $C_r<\infty$, depending on $a$ and
$\theta$, such that, for $n\ge2$,
\begin{equation}\label{eq:D-moment-bound}
   \E|D_n|^r\le C_r
   \begin{cases}
      n^{r/2},&a<1/2,\\
      (n\log n)^{r/2},&a=1/2,\\
      n^{ar},&a>1/2.
   \end{cases}
\end{equation}
For $j\ge2$, set
\begin{equation}\label{eq:mj-def}
   m_j:=\frac{aD_{j-1}}{j-1},
   \qquad m_1:=0.
\end{equation}
Then, for every $\eta>0$, there are constants $c,C>0$ such that
\begin{equation}\label{eq:localization-tail}
   \Pp\left\{\max_{\sqrt n\le j\le n}|m_j|>\eta\right\}
   \le C\exp\left\{-c n^{(1-a)/2}\right\}.
\end{equation}
When $a=0$, the probability on the left is zero.
\end{lemma}

\begin{proof}
By the martingale Rosenthal inequality
\cite[Theorem~2.12]{HallHeyde}, for every $r\ge2$,
\[
\E|\mathcal M_n|^r
\le C_r\bigg\{
\bigg(\sum_{j=1}^n\gamma_j^{-2}\bigg)^{r/2}
+\sum_{j=1}^n\gamma_j^{-r}
\bigg\}.
\]
Since $\gamma_j\ge1$,
\[
\sum_{j=1}^n\gamma_j^{-r}
\le \sum_{j=1}^n\gamma_j^{-2}
\le\bigg(\sum_{j=1}^n\gamma_j^{-2}\bigg)^{r/2}.
\]
Multiplying by $\gamma_n^r$ and using Proposition \ref{app:gamma-estimates}, we can obtain
\eqref{eq:D-moment-bound}.

Assume now that $a>0$. For \eqref{eq:localization-tail}, we apply
Freedman's inequality \cite[Theorem~1.6]{Freedman} to $\mathcal M_j-\mathcal M_1$.
Since $|\xi_j|\le1$ and $\gamma_j\ge1$, its increments are bounded by one,
and its predictable quadratic variation up to time $N$ is at most
$V_N:=\sum_{j\le N}\gamma_j^{-2}$. On a dyadic interval $N\le j\le2N$, from the relation $D_j=\gamma_j\mathcal M_j$ with $\gamma_j\le C_aN^a$, we know that the event
$|D_j|>\eta j/a$ implies $|\mathcal M_j|\ge c_\eta N^{1-a}$. Note that, since $|\mathcal M_1|\le1$, the event
$\max_{N\le j\le2N}|aD_j/j|>\eta$ implies $\max_{j\le2N}|\mathcal M_j-\mathcal M_1|\ge c_\eta N^{1-a}/2$ for large enough $N$. Moreover, by Proposition \ref{app:gamma-estimates}, we have $V_{2N}=o(N^{1-a})$ in all three regimes. Freedman's inequality, applied to $\mathcal M_j-\mathcal M_1$ and its negative with $x=c_\eta N^{1-a}/2$, thus yields
\[
\Pp\left(\max_{N\le j\le2N}|aD_j/j|>\eta\right)
\le2\exp\left\{-\frac{x^2}{2(V_{2N}+x/3)}\right\}
\le C\exp\left\{-cN^{1-a}\right\}.
\]
Note that the shift from $D_j/j$ to
$D_{j-1}/(j-1)$ in the definition of $m_j$ changes the interval endpoints
by at most one and does not affect the estimate.  By summing the resulting geometric tail over the dyadic intervals that meet $[\sqrt n,n]$, we get \eqref{eq:localization-tail} for $a>0$.
If $a=0$, then $m_j=0$ identically, so \eqref{eq:localization-tail} holds trivially.
\end{proof}

We next give the genealogical estimates needed to control local
dependence. In the copy construction, every fresh decision starts a root and
every copied decision attaches to an earlier vertex. All fresh/copy decisions
and all uniform predecessor choices are independent. Write $i\sim j$ when
$i$ and $j$ have the same root. The genealogy is independent of the Bernoulli
labels carried by the roots.

\begin{lemma}\label{lem:genealogy}
Let the root labels be Bernoulli with parameter $\beta\in(0,1)$, and set
\[
D_j^{(\beta)}:=\sum_{h\le j}(X_h-\beta),
\qquad
v_j^{(\beta)}:=\E\bigl(D_j^{(\beta)}\bigr)^2.
\]
Then
\begin{equation}\label{eq:variance-recursion}
   v_{j+1}^{(\beta)}=\left(1+\frac{2a}{j}\right)v_j^{(\beta)}+\beta(1-\beta),
\end{equation}
and for $2\le i<j$,
\begin{equation}\label{eq:root-probability-exact}
   \Pp(i\sim j)
   =\frac{a\gamma_{j-1}}{(j-1)\gamma_i}
   \bigg(1+\frac{a v_{i-1}^{(\beta)}}{\beta(1-\beta)(i-1)}\bigg).
\end{equation}
Consequently, with constants depending only on $a$,
\begin{equation}\label{eq:root-probability-bounds}
   \Pp(i\sim j)\le C
   \begin{cases}
      i^{-a}j^{a-1},&a<1/2,\\
      (\log i)i^{-1/2}j^{-1/2},&a=1/2,\\
      i^{a-1}j^{a-1},&a>1/2.
   \end{cases}
\end{equation}
Furthermore, if $A$ and $B$ are finite disjoint index sets,
$F=f((X_u)_{u\in A})$, and $G=g((X_v)_{v\in B})$, with
$|F|,|G|\le1$, then
\begin{equation}\label{eq:cov-local-functions}
   |\Cov(F,G)|\le6
   \sum_{\substack{u,v\in A\cup B\\u<v}}\Pp(u\sim v).
\end{equation}
\end{lemma}

\begin{proof}
Since
\[
   \E(X_{j+1}-\beta\mid\calF_j)=\frac{a}{j}D_j^{(\beta)},
\]
by expanding $(D_j^{(\beta)}+X_{j+1}-\beta)^2$ and taking expectation, we obtain
\[
\begin{aligned}
   v_{j+1}^{(\beta)}
   &=v_j^{(\beta)}+
     2\E\left[D_j^{(\beta)}
       \E(X_{j+1}-\beta\mid\calF_j)\right]
     +\E(X_{j+1}-\beta)^2\\
   &=\Big(1+\frac{2a}{j}\Big)v_j^{(\beta)}
     +\beta(1-\beta),
\end{aligned}
\]
where the last term uses the Bernoulli$(\beta)$ marginal law.
For $i<j$,
\[
   \E(X_j-\beta\mid\calF_i)
   =\frac{a\gamma_{j-1}}{(j-1)\gamma_i}D_i^{(\beta)}.
\]
Moreover,
\[
\begin{aligned}
   \E\bigl[(X_i-\beta)D_i^{(\beta)}\bigr]
   &=\E(X_i-\beta)^2
     +\E\left[D_{i-1}^{(\beta)}
       \E(X_i-\beta\mid\calF_{i-1})\right]\\
   &=\beta(1-\beta)+\frac{a}{i-1}v_{i-1}^{(\beta)}.
\end{aligned}
\]

On the other hand, conditional on the genealogy, $X_i$ and $X_j$ are independent unless
$i\sim j$, while on that event they have the same Bernoulli root label. Thus
\[
   \Cov(X_i,X_j)=\beta(1-\beta)\Pp(i\sim j),
\]
which proves \eqref{eq:root-probability-exact}. To bound the variance, set
$u_j=v_j^{(\beta)}/\{\beta(1-\beta)\}$. Then
\[
   u_{j+1}=\Big(1+\frac{2a}{j}\Big)u_j+1,
   \qquad u_1=1.
\]
By a simple calculation, we have
\[
   u_j=\prod_{h=1}^{j-1}\Big(1+\frac{2a}{h}\Big)
   \bigg[1+\sum_{s=1}^{j-1}
   \prod_{h=1}^{s}\Big(1+\frac{2a}{h}\Big)^{-1}\bigg].
\]
The product is comparable to $j^{2a}$, and the sum in brackets has order
$j^{1-2a}$ when $a<1/2$, order $\log j$ when $a=1/2$, and is bounded
when $a>1/2$. Hence
\[
   v_i^{(\beta)}\le C\beta(1-\beta)
   \begin{cases}
      i,&a<1/2,\\
      i\log(ei),&a=1/2,\\
      i^{2a},&a>1/2.
   \end{cases}
\]
Since $\gamma_j/\gamma_i\asymp(j/i)^a$ by Proposition \ref{app:gamma-estimates}, we prove \eqref{eq:root-probability-bounds}.

Finally, couple the root labels with independent Bernoulli$(\beta)$ labels
$(\widetilde X_u)_{u\in A\cup B}$ attached directly to the vertices, and
set $\widetilde F=f((\widetilde X_u)_{u\in A})$ and
$\widetilde G=g((\widetilde X_v)_{v\in B})$. Let $P$ be the joint law of
$(X_u)_{u\in A\cup B}$, and let $Q$ be the product Bernoulli$(\beta)$
law, that is, the joint law of $(\widetilde X_u)_{u\in A\cup B}$.
On the event that all
vertices in $A\cup B$ have distinct roots, the two label vectors coincide,
and both follow the product law. Hence, with
$d:=\sum_{u<v\in A\cup B}\Pp(u\sim v)$, the coupling characterization of
total variation gives $d_{\mathrm{TV}}(P,Q)\le d$. Since
$|\E_PH-\E_QH|\le2d$ for every $|H|\le1$, applying
this to $FG$, $F$ and $G$, and noting that $\widetilde F$ and $\widetilde G$ are independent under the product law, it follows that
\[
|\Cov(F,G)|
\le2d+|\E\widetilde F|\,|\E G-\E\widetilde G|
       +|\E G|\,|\E F-\E\widetilde F|
\le6d,
\]
which proves \eqref{eq:cov-local-functions}.
\end{proof}

\begin{lemma}\label{lem:weighted-local}
Fix $K\ge1$ and $\delta\in(0,1)$. For each $n$, let
$(B_{n,\alpha})_{\alpha\in\mathcal A_n}$ be deterministic, pairwise
disjoint subsets of $[\delta n,n]$ with $|B_{n,\alpha}|\le K$.
Let $f_{n,\alpha}$ be deterministic functions, and set
$F_{n,\alpha}=f_{n,\alpha}((X_j)_{j\in B_{n,\alpha}})$.
Suppose $|F_{n,\alpha}|\le1$. If deterministic nonnegative weights satisfy
\begin{align}\label{eq:weighted-condition}
   \sum_{\alpha}w_{n,\alpha}=O(1)
   \quad{and}\quad
   \max_\alpha w_{n,\alpha}\longrightarrow0,
\end{align}
then
\begin{equation}\label{eq:weighted-local}
   \sum_{\alpha}w_{n,\alpha}
   \bigl(F_{n,\alpha}-\E F_{n,\alpha}\bigr)\toP0.
\end{equation}
In fact, the convergence in \eqref{eq:weighted-local} holds in $L^2$.
The joint law on $B_{n,\alpha}$ converges uniformly in $\alpha$ to the
product Bernoulli$(\theta)$ law.
\end{lemma}

\begin{proof}
For indices in $[\delta n,n]$, the bound
\eqref{eq:root-probability-bounds} tends to zero uniformly. From the bound
\eqref{eq:cov-local-functions}, we therefore have
\[
   |\Cov(F_{n,\alpha},F_{n,\beta})|\le C_K r_n,
   \qquad
   r_n=\begin{cases}
      O(n^{-1}),&a<1/2,\\
      O((\log n)/n),&a=1/2,\\
      O(n^{2a-2}),&a>1/2,
   \end{cases}
\]
for $\alpha\ne\beta$. Note that the variance in \eqref{eq:weighted-local} is at most
\[
   C_K\sum_\alpha w_{n,\alpha}^2
   +C_Kr_n\left(\sum_\alpha w_{n,\alpha}\right)^2,
\]
which tends to zero by the condition \eqref{eq:weighted-condition}. The same coupling used in the proof of
Lemma~\ref{lem:genealogy} shows that the total variation distance between the
law on any one set $B_{n,\alpha}$ and the product Bernoulli law is bounded by
the probability that two of its indices share a root. This probability tends
to zero uniformly in $\alpha$. The proof is complete.
\end{proof}


\subsection{Reduction to large primes}


For $p>\sqrt n$, set $k_p=\lfloor n/p\rfloor$ and define
\begin{equation}\label{eq:prime-block}
   \calP_{n,k}:=\left\{p\in\calP:p>\sqrt n,
   \ \frac{n}{k+1}<p\le\frac nk\right\}.
\end{equation}

\begin{lemma}\label{lem:prime-block-bounds}
There is a constant $C$ such that, for $n\ge2$ and
$1\le k\le\lfloor\sqrt n\rfloor$,
\begin{align}
   \sum_{p\in\calP_{n,k}}\log p&\le C\frac nk,\label{eq:block-log}\\
   \sum_{p\in\calP_{n,k}}\log^2p&\le C\frac{n\log n}{k},\label{eq:block-log2}\\
   \sum_{p\in\calP_{n,k}}p^{\alpha-1}\log p
      &\le C_\alpha n^\alpha k^{-\alpha},
      \qquad 0<\alpha<1,\label{eq:block-alpha}\\
   \sum_{p\in\calP_{n,k}}\frac{1}{p}\log p&\le C.\label{eq:block-zero}
\end{align}
For each fixed $k$,
\begin{align}
   \frac1n\sum_{p\in\calP_{n,k}}\log p&\longrightarrow\frac1{k(k+1)},\label{eq:block-asymp-log}\\
   \frac1{n\log n}\sum_{p\in\calP_{n,k}}\log^2p&\longrightarrow\frac1{k(k+1)},\label{eq:block-asymp-log2}\\
   \frac1{n^\alpha}\sum_{p\in\calP_{n,k}}p^{\alpha-1}\log p
      &\longrightarrow\frac{k^{-\alpha}-(k+1)^{-\alpha}}{\alpha}.
      \label{eq:block-asymp-alpha}
\end{align}
\end{lemma}

\begin{proof}
The uniform bounds use only the estimates in
Proposition~\ref{app:arithmetic-estimates}. First,
\[
   \sum_{p\in\calP_{n,k}}\log p
   \le \vartheta(n/k)\le Cn/k,
\]
which proves \eqref{eq:block-log}; multiplying by $\log n$ gets
\eqref{eq:block-log2}. Since $p>n/(k+1)$ on the block and
$x^{\alpha-1}$ is decreasing,
\[
\begin{aligned}
   \sum_{p\in\calP_{n,k}}p^{\alpha-1}\log p
   &\le \left(\frac{n}{k+1}\right)^{\alpha-1}
      \sum_{p\in\calP_{n,k}}\log p\le C_\alpha n^\alpha k^{-\alpha},
\end{aligned}
\]
which is \eqref{eq:block-alpha}. The same lower bound on $p$ implies that
\[
   \sum_{p\in\calP_{n,k}}\frac{1}{p}\log p
   \le \frac{k+1}{n}\sum_{p\in\calP_{n,k}}\log p\le C,
\]
and this proves \eqref{eq:block-zero}.

Fix $k$. For all sufficiently large $n$, the condition $p>\sqrt n$ is
automatic on $(n/(k+1),n/k]$. The prime number theorem therefore yields
\[
   \frac{\vartheta(n/k)-\vartheta(n/(k+1))}{n}
   \longrightarrow \frac1k-\frac1{k+1},
\]
which is \eqref{eq:block-asymp-log}. Uniformly on this fixed-ratio interval,
$\log p/\log n\to1$, so \eqref{eq:block-asymp-log2} follows from the same
calculation. Finally, from Abel summation and the fact that $\vartheta(x)\sim x$, we have
\[
\begin{aligned}
   \sum_{n/(k+1)<p\le n/k}p^{\alpha-1}\log p
   &\sim \int_{n/(k+1)}^{n/k}x^{\alpha-1}\,\dd x=\frac{n^\alpha}{\alpha}
      \bigl(k^{-\alpha}-(k+1)^{-\alpha}\bigr),
\end{aligned}
\]
which proves \eqref{eq:block-asymp-alpha}.
\end{proof}

Let
\[
   b_n=\begin{cases}
      \sqrt{n\log n},&a\le1/2,\\
      n^a,&a>1/2\red{,}
   \end{cases}
\]
and
\[
   I_n(m):=\1_{\{A_n\cap m\N\ne\varnothing\}},
   \qquad
   S_n:=\sum_{p>\sqrt n}\log p\,I_n(p).
\]
The von Mangoldt identity gives $T_n=\sum_{m\le n}\Lambda(m)I_n(m)$.
The next deterministic estimate shows that all higher prime powers and all
primes at most $\sqrt n$ are negligible on every fluctuation scale in
Theorem~\ref{thm:fixed-main}.

\begin{lemma}\label{lem:large-prime-reduction}
We have
\begin{equation}\label{eq:large-prime-L1}
   0\le T_n-S_n
   \le \vartheta(\sqrt n)
      +\sum_{\substack{p^r\le n,\,r\ge2}}\log p
   =O\bigl(\sqrt n+n^{1/3}\log n\bigr)=o(b_n).
\end{equation}
Consequently,
\[
   T_n-S_n-\E(T_n-S_n)=o_{L^1}(b_n).
\]
\end{lemma}

\begin{proof}
In the von Mangoldt expansion of $T_n$, the sum $S_n$ retains only the first
powers of primes larger than $\sqrt n$. Every omitted term is nonnegative.
The omitted first powers contribute at most $\vartheta(\sqrt n)$, and all
omitted powers of exponent at least two contribute at most $\sum_{\substack{p^r\le n,\,r\ge2}}\log p$.
By Proposition~\ref{app:arithmetic-estimates}, the first quantity is
$O(\sqrt n)$ and the second quantity is $O(\sqrt n+n^{1/3}\log n)$, which are $o(b_n)$ in all three memory regimes.
This proves \eqref{eq:large-prime-L1}. The centered $L^1$ assertion is obvious.
\end{proof}

The blockwise law of large numbers and the reduction to large primes now
identify the first-order limit of $T_n$.

\begin{proposition}\label{prop:fixed-lln}
For fixed $\theta\in(0,1)$ and $0\le a<1$, ${T_n}/{n}\toP c(\theta)$.
\end{proposition}

\begin{proof}
From Lemma~\ref{lem:large-prime-reduction}, it is enough to consider $S_n$. Fix $K\geq1$. If $p\in\calP_{n,k}$ with
$k\le K$, then $I_n(p)$ depends on
$(X_p,X_{2p},\ldots,X_{kp})$. These supports are pairwise disjoint as $p$
varies because an integer at most $n$ cannot have two prime divisors larger
than $\sqrt n$. By Lemma~\ref{lem:weighted-local} and
\eqref{eq:block-asymp-log}, we have
\[
   \frac1n\sum_{k=1}^K\sum_{p\in\calP_{n,k}}\log p\,I_n(p)
   \toP
   \sum_{k=1}^K\frac{1-q^k}{k(k+1)}.
\]
In fact, Lemma~\ref{lem:weighted-local} also implies
\[
   \sup_{p\in\calP_{n,k}}
   \left|\E I_n(p)-(1-q^k)\right|\longrightarrow0
\]
for each fixed $k$.
The omitted primes satisfy $p\le n/(K+1)$, so their total contribution is at
most $\vartheta(n/(K+1))\le Cn/K$. Hence, letting $K\to\infty$ we obtain
\[
   \frac{S_n}{n}\toP\sum_{k\ge1}\frac{1-q^k}{k(k+1)}.
\]
Finally, since $\sum_{k\ge1}q^k/k=-\log(1-q)$, it follows that
\[
   \sum_{k\ge1}\frac{1-q^k}{k(k+1)}
   =\frac{1-q}{q}\log\frac1{1-q}
   =\frac{\theta\log(1/\theta)}{1-\theta}.
\]
The proof is complete.
\end{proof}

We now localize the memory drift and derive the exact decomposition of the
large-prime contribution. Recall that $\xi_j=X_j-\E(X_j\mid\calF_{j-1})$ and $m_j=aD_{j-1}/(j-1)$, so that
$X_j-\theta=m_j+\xi_j$. Fix a small $\eta>0$ and set
\[
   \tau_n:=\inf\{j\ge\lceil\sqrt n\rceil:\; |m_j|>\eta\},
   \qquad \inf\varnothing=\infty.
\]
Let $G_n=\{\tau_n>n\}$ and $\chi_j=\1_{\{\tau_n>j\}}$. The variable
$\chi_j$ is $\calF_{j-1}$-measurable. Define
\[
   \widehat m_j:=\chi_jm_j,
   \qquad
   \widehat \xi_j:=\chi_j\xi_j,
   \qquad
   \widehat X_j:=\theta+\widehat m_j+\widehat \xi_j.
\]
Before $\tau_n$, the modified row equals the original row. At and after
$\tau_n$, it equals the constant $\theta$. Hence $0\le\widehat X_j\le1$.
Define $\widehat S_n$ from the modified row by replacing $I_n(p)$ with
\[
   \widehat I_n(p):=1-\prod_{r=1}^{k_p}(1-\widehat X_{rp}).
\]
On $G_n$, the original and modified rows coincide up to time $n$. Since
$S_n$ and $\widehat S_n$ are bounded by $\vartheta(n)=O(n)$,
by Lemma~\ref{lem:count-moments-localization} we have, for $r=1,2$,
\begin{equation}\label{eq:localization-transfer}
   \|S_n-\widehat S_n\|_r
   \le Cn\,\Pp(G_n^c)^{1/r}=o(b_n),
\end{equation}
and in particular
$\|(S_n-\E S_n)-(\widehat S_n-\E\widehat S_n)\|_1=o(b_n)$.
We retain the hats throughout the localized argument. The
symbol $D_j$ always denotes the original count process, and identities
between $\widehat m_j$ and $D_{j-1}$ are used only on $G_n$. For each fixed $K$, the finite-row sums and products used below in the first $K$ blocks are bounded by $C_Kn^{C_K}$ pathwise, while $\Pp(G_n^c)$ decays
exponentially by Lemma~\ref{lem:count-moments-localization}. H\"older's inequality therefore
yields $\|R_n\1_{G_n^c}\|_r\le C_Kn^{C_K}\Pp(G_n^c)^{1/r}$ for every such $R_n$ and every $r\ge1$, and
$\E(P_n|W_{a,\theta}|\1_{G_n^c})\le P_n\|W_{a,\theta}\|_2\Pp(G_n^c)^{1/2}$ for every deterministic
nonnegative polynomial $P_n$. Both bounds decay faster than any power of $n$.

Choose $\rho\in(q,1)$. By taking $\eta$ small enough, for $r=2,4$ and
$j\ge\lceil\sqrt n\rceil$, we can get
\begin{equation}\label{eq:product-contraction-one}
   \E\bigl(|q-\widehat\xi_j|^r\mid\calF_{j-1}\bigr)\le\rho
   \quad\text{and}\quad
   \E\bigl(|1-\widehat X_j|^r\mid\calF_{j-1}\bigr)\le\rho.
\end{equation}
Indeed, the conditional success probability is $p_j=\theta+m_j$. On this range and before the stopping time, the variable
$q-\widehat\xi_j$ equals $m_j$ when $X_j=1$ and equals $1+m_j$ when
$X_j=0$. Hence
\[
   \E\bigl(|q-\widehat\xi_j|^r\mid\calF_{j-1}\bigr)
   =p_j|m_j|^r+(1-p_j)|1+m_j|^r.
\]
As a function of $m_j$, the right-hand side is continuous and equals $q$ at zero. Similarly,
\[
   \E\bigl(|1-\widehat X_j|^r\mid\calF_{j-1}\bigr)=q-m_j\le q+\eta.
\]
After the stopping time, the two quantities are $q^r$. Thus we can choose $\eta$ small enough so that both conditional moments are bounded by $\rho$. If factors at distinct times are arranged in increasing order, successive conditioning
then supplies one factor $\rho$ for every factor raised to the second or
fourth power in their product. Every time used below has the form $rp$ with
$p>\sqrt n$, so the restriction $j\ge\lceil\sqrt n\rceil$ covers every
application of \eqref{eq:product-contraction-one}.

For $p\in\calP_{n,k}$, write
\[
   \widehat A_{r,p}:=q-\widehat\xi_{rp}
   \quad\text{and}\quad
   \widehat B_{r,p}:=1-\widehat X_{rp}
   =\widehat A_{r,p}-\widehat m_{rp}.
\]
The exact telescoping identities
\begin{align}
   \widehat H_{n,p}&:=q^k-\prod_{r=1}^k\widehat A_{r,p}
   =\sum_{\ell=1}^k\widehat\xi_{\ell p}q^{k-\ell}
      \prod_{r<\ell}\widehat A_{r,p},\label{eq:H-exact}\\
   \widehat V_{n,p}&:=\prod_{r=1}^k\widehat A_{r,p}
      -\prod_{r=1}^k\widehat B_{r,p}=\sum_{\ell=1}^k\widehat m_{\ell p}
      \prod_{r<\ell}\widehat B_{r,p}
      \prod_{r>\ell}\widehat A_{r,p}
      \label{eq:V-exact}
\end{align}
show that
\begin{equation}\label{eq:S-exact-decomposition}
   \widehat S_n-A_n^{(0)}=\widehat U_n+\widehat V_n,
\end{equation}
where
\begin{align*}
   A_n^{(0)}:=\sum_{p>\sqrt n}(1-q^{k_p})\log p,\quad
   \widehat U_n:=\sum_{p>\sqrt n}\widehat H_{n,p}\log p,\quad
   \widehat V_n:=\sum_{p>\sqrt n}\widehat V_{n,p}\log p.
\end{align*}
Note that an integer $t\leq n$ has at most one representation $t=\ell p$ with $p>\sqrt n$. Hence, grouping \eqref{eq:H-exact} by $t=\ell p$, we have $\widehat U_n=\sum_{t\le n}c_{n,t}\widehat\xi_t$, where
\[
   c_{n,t}:=
   \begin{cases}
      (\log p)q^{k_p-\ell}\displaystyle\prod_{r<\ell}\widehat A_{r,p},
         &t=\ell p,\quad p>\sqrt n,\quad 1\le\ell\le k_p,\\[1mm]
      0,&\text{if no such representation exists}.
   \end{cases}
\]
Moreover, every factor in $c_{n,t}$ is observed before time $t$, so $c_{n,t}$ is
$\calF_{t-1}$-measurable. Since $\widehat\xi_t$ is a martingale
difference, the partial sums
\begin{align}\label{eq:martingale-hat}
  \widehat U_{n,t}:=\sum_{s\le t}c_{n,s}\widehat\xi_s,\qquad 0\le t\le n,
\end{align}
then form a martingale array, and $\widehat U_{n,n}=\widehat U_n$.
The term $\widehat V_n$ contains the contribution of the predictable
drifts $\widehat m_j$, whose coefficients may involve later innovations, unlike those in $\widehat U_n$.
We therefore treat it separately, by the count projection in the next subsection.


\subsection{Martingale limits and count projection}


For square-integrable martingales $R$ and $S$ adapted to $(\calF_t)$,
write $\Delta R_t=R_t-R_{t-1}$ and
\[
   \langle R,S\rangle_n
      :=\sum_{t=1}^n\E(\Delta R_t\Delta S_t\mid\calF_{t-1}),
   \qquad \langle R\rangle_n:=\langle R,R\rangle_n.
\]
These are the predictable covariation and quadratic variation.

We first identify the fluctuation carried by the innovation martingale. The
same calculation also determines its joint behavior with the retained count
at the critical value.

\begin{proposition}\label{prop:innovation-limit}
If $0\le a\le1/2$, then
\begin{equation}\label{eq:U-clt}
   \frac{\widehat U_n}{\sqrt{n\log n}}\tod N(0,g(q)).
\end{equation}
At $a=1/2$,
\begin{equation}\label{eq:joint-U-D}
   \left(\frac{\widehat U_n}{\sqrt{n\log n}},
   \frac{D_n}{\sqrt{n\log n}}\right)
   \tod (G,\sqrt{\theta q}\,Z),
\end{equation}
where $G\sim N(0,g(q))$, $Z\sim N(0,1)$, and $G$ and $Z$ are
independent. If $1/2<a<1$, then $\widehat U_n/n^a\toP0$.
\end{proposition}

\begin{proof}
For $K\ge1$, let $\widehat U_{n,\le K}$ be the contribution of the
blocks $k\le K$, and set
$\widehat U_{n,>K}:=\widehat U_n-\widehat U_{n,\le K}$. For the martingale array $\widehat U_{n,t}$ defined in \eqref{eq:martingale-hat}, by restricting its coefficients to the first $K$ blocks, we have
\begin{align}
   \widehat U_{n,t}^{(K)}
   :=\sum_{s\le t}c_{n,s}\1_{\{k_p\le K\}}\,\widehat\xi_s,
   \qquad 0\le t\le n,
\end{align}
which is again an $\calF_t$-martingale starting at zero, with terminal value
$\widehat U_{n,n}^{(K)}=\widehat U_{n,\le K}$.

To bound $\widehat U_{n,>K}$, expand each $\widehat H_{n,p}$ by the
telescoping identity \eqref{eq:H-exact}. From the orthogonality of the martingale
differences, \eqref{eq:product-contraction-one} and \eqref{eq:block-log2}, we then obtain
\begin{equation}\label{eq:U-tail-L2}
\begin{aligned}
   \E|\widehat U_{n,>K}|^2
   &=\sum_{k>K}\sum_{p\in\calP_{n,k}}\log^2p
     \sum_{\ell=1}^k q^{2(k-\ell)}
     \E\Bigl(\widehat\xi_{\ell p}^2\prod_{r<\ell}\widehat A_{r,p}^2\Bigr)\\
   &\le\sum_{k>K}\sum_{p\in\calP_{n,k}}\log^2p
     \sum_{\ell=1}^k q^{2(k-\ell)}\rho^{\ell-1}\\
     &\le \sum_{k>K}k\rho^{k-1}
       \sum_{p\in\calP_{n,k}}\log^2p\leq Cn\log n\sum_{k>K}\rho^{k-1}.
\end{aligned}
\end{equation}
Here, we use the fact $q^2<\rho$ in the second inequality. The same calculation without the restriction $k>K$ can give
$\E\widehat U_n^2=O(n\log n)$, so $\widehat U_n/n^a\toP0$ when $a>1/2$.

For fixed $K$, the predictable quadratic variation of the
truncated array $(\widehat U_{n,t}^{(K)})_{0\le t\le n}$ is
\begin{align}\label{eq:QV-hat}
   \langle\widehat U_{n,\cdot}^{(K)}\rangle_n=\sum_{k=1}^K\sum_{p\in\calP_{n,k}}\log^2p
   \sum_{\ell=1}^kq^{2(k-\ell)}
   \E(\widehat\xi_{\ell p}^2\mid\calF_{\ell p-1})
   \prod_{r<\ell}\widehat A_{r,p}^2.
\end{align}
Note that all indices $\ell p$ are at least $n/(K+1)$ and all products have a deterministic bound. By \eqref{eq:localization-tail},
\begin{align}\label{eq:martingale-max}
   M_{n,K}^*:=\max_{n/(K+1)\le j\le n}|m_j|\toP0.
\end{align}
On $G_n$, by the telescoping inequality for finite products and a simple calculation, we obtain
\[
 \Big|\prod_{r<\ell}\widehat A_{r,p}^2
       -\prod_{r<\ell}(1-X_{rp})\Big|
 \le C_KM_{n,K}^*
\]
and
\[
 \sup_{n/(K+1)\le j\le n}
 \left|\E(\widehat\xi_j^2\mid\calF_{j-1})-\theta q\right|
 \le CM_{n,K}^*.
\]

Define
\[
   Q^{(U)}_{n,K}:=\theta q
   \sum_{k=1}^K\sum_{p\in\calP_{n,k}}\log^2p
   \sum_{\ell=1}^kq^{2(k-\ell)}
      \prod_{r<\ell}(1-X_{rp}).
\]
The two replacements above are valid on $G_n$, with error
$O_K(M_{n,K}^*)$ per summand. On $G_n^c$, the products are bounded
deterministically, so both
$\langle\widehat U_{n,\cdot}^{(K)}\rangle_n$ and $Q^{(U)}_{n,K}$
are at most the total $\log^2p$ weight, namely $O_K(n\log n)$.  Hence, by \eqref{eq:localization-tail} and \eqref{eq:martingale-max}, we have
\[
   \frac{\bigl|
      \langle\widehat U_{n,\cdot}^{(K)}\rangle_n-Q^{(U)}_{n,K}
   \bigr|}{n\log n}
   \le C_K\bigl(\1_{G_n^c}+M_{n,K}^*\bigr)\toP0.
\]
It suffices to identify the limit of $Q^{(U)}_{n,K}$. For fixed
$p\in\calP_{n,k}$, define
\[
   h^{(U)}_{n,p}(x_p,x_{2p},\ldots,x_{kp})
   :=\sum_{\ell=1}^kq^{2(k-\ell)}\prod_{r<\ell}(1-x_{rp}), \quad x_{\ell p}\in[0,1],\; 1\le\ell\le k.
\]
Note that the function $h^{(U)}_{n,p}$ is bounded by $(1-q^2)^{-1}$, and the sets $\{p,2p,\ldots,kp\}$  are disjoint for distinct primes. So applying Lemma~\ref{lem:weighted-local} with weights
$w_{n,p}=\log^2p/(n\log n)$, we get
\[
\frac{\bigl( Q^{(U)}_{n,K}-\E Q^{(U)}_{n,K}\bigr)}{n\log n}\toP0.
\]
Under the product Bernoulli$(\theta)$ law,
\[
   \E h^{(U)}_{n,p}(X_p,X_{2p},\ldots,X_{kp})
   =\sum_{\ell=1}^kq^{2(k-\ell)}q^{\ell-1},
\]
so, by the last assertion of Lemma~\ref{lem:weighted-local}, the
expectation of the bounded function $h^{(U)}_{n,p}$ under the original law
converges to this product-law expectation, uniformly over the first $K$
blocks. It follows from \eqref{eq:block-asymp-log2} that
\[
   \frac{Q^{(U)}_{n,K}}{n\log n}
   \toP
   g_K(q):=\sum_{k=1}^K\frac{q^k(1-q^k)}{k(k+1)}.
\]
For fixed $K$, each increment of the martingale $\widehat U_{n, \le K}$
is bounded by $C_K\log n$, which is $o(\sqrt{n\log n})$, so the conditional Lindeberg condition holds. The martingale central limit theorem \cite[Corollary~3.1]{HallHeyde} then yields that $\widehat U_{n, \le K}/\sqrt{n\log n}\tod N(0,g_K(q))$. Moreover, for every fixed $K$, by Chebyshev's inequality and \eqref{eq:U-tail-L2}, we have
\[
   \Pp\bigl(|\widehat U_{n,>K}|>\varepsilon\sqrt{n\log n}\bigr)
   \le C\varepsilon^{-2}\rho^K,
\]
so applying the converging-together theorem \cite[Theorem~3.2]{Billingsley} and noting
that $g_K(q)\to g(q)$, we prove \eqref{eq:U-clt}.

At $a=1/2$, define the martingale array
\[
   \widetilde D_{n,t}:=\gamma_n\mathcal M_t=\gamma_n\sum_{j=1}^t
   \frac{\xi_j}{\gamma_j},
   \qquad 0\le t\le n,
\]
with $\widetilde D_{n,0}=0$ and $\widetilde D_{n,n}=D_n$ by the representation $D_j=\gamma_j\mathcal M_j$. Since $|\xi_j|\le1$, $\gamma_j\ge1$, and
$\gamma_n\le C\sqrt n$ by Proposition~\ref{app:gamma-estimates}, its largest increment after division by
$\sqrt{n\log n}$ is $O((\log n)^{-1/2})$. The predictable quadratic variation of
$\widetilde D_{n,\cdot}$ at time $n$ satisfies
\begin{equation}\label{eq:D-bracket}
   \frac{\langle\widetilde D_{n,\cdot}\rangle_n}{n\log n}
   =\frac{\gamma_n^2}{n\log n}
    \sum_{j=1}^n\frac{p_j(1-p_j)}{\gamma_j^2}
   \toP\theta q.
\end{equation}
To see this, setting
$w_{nj}:=\gamma_j^{-2}\big/\sum_{i\le n}\gamma_i^{-2}$ and recalling that $p_j=\theta+m_j$, we split the
sum as follows,
\[
   \frac{\langle\widetilde D_{n,\cdot}\rangle_n}{n\log n}
   =\theta q\,\frac{\gamma_n^2}{n\log n}
      \sum_{j\le n}\gamma_j^{-2}
   +\frac{\gamma_n^2}{n\log n}
    \sum_{i\le n}\gamma_i^{-2}
    \sum_{j\le n}w_{nj}\bigl((q-\theta)m_j-m_j^2\bigr).
\]
The first term converges to $\theta q$, since
$\gamma_n^2/n\to\Gamma(3/2)^{-2}$ and
$\sum_{j\le n}\gamma_j^{-2}\sim\Gamma(3/2)^2\log n$ by
Proposition~\ref{app:gamma-estimates}. In the second term, the
prefactor is $O(1)$, and the sequence $(q-\theta)m_j-m_j^2$ is bounded
and tends to zero almost surely because $D_j/j\to0$ almost surely. Moreover,
the weights $w_{nj}$ are nonnegative with total mass one, and
$\max_{j\le n}w_{nj}=O(1/\log n)$. By the Toeplitz lemma the second term therefore tends
to zero almost surely. This proves \eqref{eq:D-bracket}.

We now compute the predictable cross variation between
$\widehat U_{n,\le K}$ and $D_n$.
For $k\le K$ and $p\in\calP_{n,k}$, every time $\ell p$ is at least
$n/(K+1)$. Hence $\gamma_n/\gamma_{\ell p}\le C_K$, the products in
\eqref{eq:H-exact} are bounded by a constant depending only on $K$, and
\[
   |\langle\widehat U_{n,\cdot}^{(K)},
      \widetilde D_{n,\cdot}\rangle_n|
   \le C_K\sum_{k\le K}\sum_{p\in\calP_{n,k}}\log p
   =O_K(n).
\]
After division by $n\log n$, the cross variation tends to zero. Note that the largest
normalized increment of each coordinate tends to zero, and hence the
vector conditional Lindeberg condition holds. Applying the Cram\'er--Wold device and the
martingale central limit theorem \cite[Corollary~3.1]{HallHeyde}, we then obtain a
centered Gaussian limit with diagonal covariance matrix
$\operatorname{diag}(g_K(q),\theta q)$. Finally, letting $K\to\infty$ and using \eqref{eq:U-tail-L2}, by the
converging-together theorem again, we prove \eqref{eq:joint-U-D} and the independence of its two limiting coordinates.
\end{proof}

It remains to identify the memory term. To project the drift onto the terminal count, we need uniform approximations
of $D_j$ over time windows whose indices are proportional to $n$.

\begin{lemma}\label{lem:count-window}
At $a=1/2$, for every $\delta\in(0,1)$,
\begin{equation}\label{eq:critical-window}
   \sup_{\delta n\le j\le n}
   \frac{|D_j-(j/n)^{1/2}D_n|}{\sqrt{n\log n}}
   \longrightarrow0
   \quad\text{in }L^2.
\end{equation}
If $1/2<a<1$, then
\begin{equation}\label{eq:strong-window}
   \sup_{\delta n\le j\le n}
   \Big|\frac{D_j}{j^a}-W_{a,\theta}\Big|
   \longrightarrow0
   \quad\text{in }L^2.
\end{equation}
\end{lemma}

\begin{proof}
For $a=1/2$, write $D_j=\gamma_j\mathcal M_j$ and set
$j_0=\lfloor\delta n\rfloor$. Since
\[
   \sup_{j_0\le j\le n}|\mathcal M_n-\mathcal M_j|
   \le |\mathcal M_n-\mathcal M_{j_0}|
      +\sup_{j_0\le j\le n}|\mathcal M_j-\mathcal M_{j_0}|,
\]
by Doob's inequality \cite[Theorem~2.1]{HallHeyde}, we have
\[
   \E\sup_{j_0\le j\le n}|\mathcal M_n-\mathcal M_j|^2
   \le C\sum_{j_0<h\le n}\gamma_h^{-2}=O_\delta(1).
\]
Note that, for $\delta n\le j\le n$,
\[
   D_j-\Big(\frac jn\Big)^{1/2}D_n
   =\gamma_j(\mathcal M_j-\mathcal M_n)
    +\Big\{\frac{\gamma_j}{\gamma_n}
      -\Big(\frac jn\Big)^{1/2}\Big\}D_n.
\]
Then, from Proposition \ref{prop:count-transition} and Proposition \ref{app:gamma-estimates}, we obtain \eqref{eq:critical-window}.

If $1/2<a<1$, the martingale $\mathcal M_j$ is $L^2$ bounded and hence converges in $L^2$ to $\mathcal M_\infty$.
More explicitly, for every integer $N$,
\[
   \Big\|\sup_{j\ge N}|\mathcal M_j-\mathcal M_\infty|\Big\|_2
   \le 3\|\mathcal M_N-\mathcal M_\infty\|_2\longrightarrow0.
\]
Furthermore,
\begin{align*}
\sup_{\delta n\le j\le n}
   \Big|\frac{D_j}{j^a}-W_{a,\theta}\Big|\le\sup_{\delta n\le j\le n}|\mathcal M_j-\mathcal M_n |\frac{\gamma_j}{j^a}
   +|\mathcal M_n|\sup_{\delta n\le j\le n}\Big|\frac{\gamma_j}{j^a}-\frac{\gamma_n}{n^a}\Big|+\Big|\frac{D_n}{n^a}-W_{a,\theta}\Big|.
\end{align*}
Since $\gamma_j/j^a\to1/\Gamma(1+a)$ uniformly on $\delta n\le j\le n$, the first term on the right-hand side tends to zero in $L^2$. The second term also tends to zero in $L^2$ by the uniform gamma-ratio
estimate in Proposition~\ref{app:gamma-estimates}, while the third term
tends to zero in $L^2$ by Proposition~\ref{prop:count-transition}.
This proves \eqref{eq:strong-window}.
\end{proof}

The window estimates allow the local memory drifts to be replaced by the
terminal retained count. The following proposition carries out this
projection and identifies the coefficient in each memory regime.

\begin{proposition}\label{prop:drift-projection}
The drift term in \eqref{eq:S-exact-decomposition} satisfies
\begin{align}
0\le a<1/2:
&\qquad \frac{\widehat V_n}{\sqrt{n\log n}}\longrightarrow0
   \quad\text{in }L^1,\label{eq:V-weak}\\
 a=1/2:
&\qquad \frac{\widehat V_n-K_{1/2,\theta}D_n}{\sqrt{n\log n}}
   \longrightarrow0
   \quad\text{in }L^1,\label{eq:V-critical}\\
1/2<a<1:
&\qquad \frac{\widehat V_n}{n^a}\longrightarrow K_{a,\theta}W_{a,\theta}
   \quad\text{in }L^1.\label{eq:V-strong}
\end{align}
\end{proposition}

\begin{proof}
For $p\in\calP_{n,k}$ and $1\le\ell\le k$, let
\begin{equation}\label{eq:C-coefficient}
   \widehat C_{n,p,\ell}:=
   \prod_{r<\ell}\widehat B_{r,p}
   \prod_{r>\ell}\widehat A_{r,p}.
\end{equation}
We apply \eqref{eq:product-contraction-one} successively to the
factors in \eqref{eq:C-coefficient} and obtain
\begin{equation}\label{eq:C-moment}
   \|\widehat C_{n,p,\ell}\|_2\le\rho^{(k-1)/2}.
\end{equation}
For $0\le a<1/2$, \eqref{eq:D-moment-bound} implies
\[
   \|\widehat m_{\ell p}\|_2
   =\frac{a\|\chi_{\ell p}D_{\ell p-1}\|_2}{\ell p-1}
   \le C(\ell p)^{-1/2}.
\]
Combining this with \eqref{eq:C-moment} and the Cauchy--Schwarz
inequality, we have
\[
\begin{aligned}
   \E|\widehat V_n|
   &\le
   \sum_{k\ge1}\sum_{p\in\calP_{n,k}}\log p
   \sum_{\ell=1}^k
   \E|\widehat m_{\ell p}\widehat C_{n,p,\ell}|\\
   &\le
   C\sum_{k\ge1}\rho^{(k-1)/2}
   \Big(\sum_{p\in\calP_{n,k}}p^{-1/2}\log p\Big)
   \sum_{\ell=1}^k\ell^{-1/2}\\
   &\le
   C\sum_{k\ge1}\rho^{(k-1)/2}
   \bigl(\sqrt n\,k^{-1/2}\bigr)\sqrt k
   =O(\sqrt n),
\end{aligned}
\]
where we used \eqref{eq:block-alpha} with $\alpha=1/2$.
This proves \eqref{eq:V-weak}.

For the critical and strong regimes, fix $K\ge1$, let
$\widehat V_{n,\le K}$ denote the contribution from the blocks
$k\le K$, and set
$\widehat V_{n,>K}=\widehat V_n-\widehat V_{n,\le K}$. By
\eqref{eq:D-moment-bound} again,
\[
   \|\widehat m_{\ell p}\|_2
   \le
   \begin{cases}
      C\sqrt{\log n}\,(\ell p)^{-1/2},&a=\dfrac12,\\[4pt]
      C(\ell p)^{a-1},&\dfrac12<a<1.
   \end{cases}
\]
The same Cauchy--Schwarz argument, combined with
\eqref{eq:C-moment} and \eqref{eq:block-alpha}, now shows that
\begin{equation}\label{eq:V-tail}
   \|\widehat V_{n,>K}\|_1
   \le Cb_n\sum_{k>K}\rho^{(k-1)/2},
\end{equation}
where $C$ is independent of $n$ and $K$. Thus it suffices to analyze $\widehat V_{n,\le K}$ for fixed $K$,
and then let $K\to\infty$.

For $a\ge1/2$ and $K\ge1$, define
\begin{equation}\label{eq:K-truncated}
   K_{a,\theta}^{(K)}:=
   \sum_{k=1}^Kq^{k-1}\bigl(k^{-a}-(k+1)^{-a}\bigr)
   \sum_{\ell=1}^k\ell^{a-1},
\end{equation}
and
\begin{equation}\label{eq:random-drift-coeff}
   Q_{n,K}^{(a)}:=\frac{a}{n^a}
   \sum_{k=1}^K\sum_{p\in\calP_{n,k}}
   p^{a-1}\log p
   \sum_{\ell=1}^k\ell^{a-1}\widehat C_{n,p,\ell}.
\end{equation}
We claim that
\begin{equation}\label{eq:random-coeff-limit}
   Q_{n,K}^{(a)}\longrightarrow K_{a,\theta}^{(K)}
   \quad\text{in }L^2.
\end{equation}

To prove \eqref{eq:random-coeff-limit}, fix $K\ge1$. All relevant indices are proportional to $n$. Since $D_j/j\to0$ almost surely, we have $\max_{j\ge n/(K+1)}|m_j|\to0$ almost surely. Set
\[
   M_{n,K}^*:=\max_{n/(K+1)\le j\le n}|m_j|
\]
and define the approximating coefficient
\[
   \widetilde Q_{n,K}^{(a)}:=\frac{a}{n^a}
   \sum_{k=1}^K\sum_{p\in\calP_{n,k}}
   p^{a-1}\log p\sum_{\ell=1}^k\ell^{a-1}
      \prod_{r\ne\ell}(1-X_{rp}).
\]
Since $p_j=\theta+m_j\in[0,1]$, we have
$-\theta\le m_j\le1-\theta$. Thus the variables $M_{n,K}^*$ are uniformly
bounded and converge to zero almost surely. On $G_n$, the telescoping inequality for finite products implies
\[
   \Big|\widehat C_{n,p,\ell}
      -\prod_{r\ne\ell}(1-X_{rp})\Big|
   \le C_KM_{n,K}^*.
\]
By \eqref{eq:block-alpha}, the normalized weights satisfy
\[
   \frac1{n^a}\sum_{p\in\calP_{n,k}}
      p^{a-1}\log p=O_K(1).
\]
Bounding the contribution on $G_n^c$ by an indicator, we therefore
obtain
\begin{equation}\label{eq:Q-approximation}
   |Q_{n,K}^{(a)}-\widetilde Q_{n,K}^{(a)}|
   \le C_K\bigl(M_{n,K}^*+\1_{G_n^c}\bigr)
   \longrightarrow0
   \quad\text{in }L^2.
\end{equation}
Indeed, $M_{n,K}^*$ is uniformly bounded and tends to zero almost
surely, while $\Pp(G_n^c)\to0$.

For fixed $p\in\calP_{n,k}$, define
\[
   h^{(a)}_{n,p}(x_p,x_{2p},\ldots,x_{kp})
   :=\sum_{\ell=1}^k\ell^{a-1}
      \prod_{r\ne\ell}(1-x_{rp}), \quad x_{\ell p}\in[0,1],\; 1\le\ell\le k.
\]
Each $h^{(a)}_{n,p}$ is bounded by a constant depending only on $a$ and
$K$, and the sets $\{p,2p,\ldots,kp\}$ are disjoint for distinct primes. After dividing $h^{(a)}_{n,p}$ by its uniform bound, we apply
Lemma~\ref{lem:weighted-local} separately for each $k\le K$, with
weights $w_{n,p}=n^{-a}p^{a-1}\log p$, whose sum is bounded and
whose maximum tends to zero. Note that, under the product
Bernoulli$(\theta)$ law,
\[
   \E h^{(a)}_{n,p}(X_p,X_{2p},\ldots,X_{kp})
   =q^{k-1}\sum_{\ell=1}^k\ell^{a-1}.
\]
Therefore, by Lemma~\ref{lem:weighted-local} and
\eqref{eq:block-asymp-alpha},
\[
\begin{aligned}
   \widetilde Q_{n,K}^{(a)}
   &\longrightarrow
   a\sum_{k=1}^Kq^{k-1}\sum_{\ell=1}^k\ell^{a-1}
   \frac{k^{-a}-(k+1)^{-a}}{a}=K_{a,\theta}^{(K)}
   \quad\text{in }L^2.
\end{aligned}
\]
Combining this with \eqref{eq:Q-approximation}, we prove \eqref{eq:random-coeff-limit}.

We now treat the critical and strong regimes together. Define
\[
   \zeta_{n,a}:=
   \begin{cases}
      D_n/n^a,&a=\dfrac12,\\[4pt]
      W_{a,\theta},&\dfrac12<a<1.
   \end{cases}
\]
With this notation, the desired estimates \eqref{eq:V-critical} and
\eqref{eq:V-strong} both take the form
\[
   \widehat V_n=n^aK_{a,\theta}\zeta_{n,a}+o_{L^1}(b_n).
\]

Fix $K\ge1$. For all indices in the first $K$ blocks, define
\[
   R_{n,p,\ell}:=\widehat m_{\ell p}
      -a\zeta_{n,a}(\ell p)^{a-1}.
\]
Let $h=\ell p-1$ and $\delta_K=1/\{2(K+1)\}$. For sufficiently
large $n$, all relevant indices satisfy $h\ge\delta_Kn$. On $G_n$,
we have
\[
\begin{aligned}
   R_{n,p,\ell}
   &=ah^{a-1}\Big(\frac{D_h}{h^a}-\zeta_{n,a}\Big)
     +a\zeta_{n,a}
       \left\{h^{a-1}-(h+1)^{a-1}\right\}.
\end{aligned}
\]
By the mean value theorem,
\[
\begin{aligned}
   \max_{\substack{k\le K,\ p\in\calP_{n,k}\\1\le\ell\le k}}
   |R_{n,p,\ell}|
   &\le C_Kn^{a-1}
      \sup_{\delta_Kn\le h\le n}
      \Big|\frac{D_h}{h^a}-\zeta_{n,a}\Big|
      +C_Kn^{a-2}|\zeta_{n,a}|.
\end{aligned}
\]
Lemma~\ref{lem:count-window}, together with
\eqref{eq:D-moment-bound} at $a=1/2$ and
$W_{a,\theta}\in L^2$ when $a>1/2$, yields
\[
   \Big\|
      \sup_{\delta_Kn\le h\le n}
      \Big|\frac{D_h}{h^a}-\zeta_{n,a}\Big|
   \Big\|_2
   =o\!\left(\frac{b_n}{n^a}\right),
   \qquad
   \|\zeta_{n,a}\|_2=O\!\left(\frac{b_n}{n^a}\right).
\]
Consequently, uniformly over the first $K$ blocks,
\begin{equation}\label{eq:R-good-event}
   \max_{\substack{k\le K,\ p\in\calP_{n,k}\\1\le\ell\le k}}
   \|R_{n,p,\ell}\1_{G_n}\|_2
   =o\!\left(\frac{b_n}{n}\right).
\end{equation}

On $G_n^c$, the uniform boundedness of $\widehat m_{\ell p}$ implies that
\[
\begin{aligned}
   \|R_{n,p,\ell}\1_{G_n^c}\|_2
   &\le C\Pp(G_n^c)^{1/2}
      +C_Kn^{a-1}
       \|\zeta_{n,a}\1_{G_n^c}\|_2,
\end{aligned}
\]
uniformly over the same indices. The first term is
super-polynomially small by \eqref{eq:localization-tail}. At
$a=1/2$, H\"older's inequality, \eqref{eq:D-moment-bound} with
$r=4$, and the localization tail show that the second term is
$o(b_n/n)$. For $a>1/2$, the same conclusion follows from dominated
convergence applied to $W_{a,\theta}\1_{G_n^c}$. Combining these
bounds with \eqref{eq:R-good-event}, we obtain
\begin{equation}\label{eq:R-uniform}
   \max_{\substack{k\le K,\ p\in\calP_{n,k}\\1\le\ell\le k}}
   \|R_{n,p,\ell}\|_2
   =o\!\left(\frac{b_n}{n}\right).
\end{equation}

We now substitute the approximation
\[
   \widehat m_{\ell p}
   =a\zeta_{n,a}(\ell p)^{a-1}+R_{n,p,\ell}
\]
into $\widehat V_{n,\le K}$. By the Cauchy--Schwarz inequality,
\eqref{eq:C-moment}, \eqref{eq:R-uniform}, and
$\sum_{k\le K}\sum_{p\in\calP_{n,k}}\log p=O_K(n)$, the total
contribution of the residuals is bounded in $L^1$ by
\[
   C_Kn
   \max_{\substack{k\le K,\ p\in\calP_{n,k}\\1\le\ell\le k}}
   \|R_{n,p,\ell}\|_2
   =o(b_n).
\]
It follows that, for each fixed $K$,
\begin{equation}\label{eq:V-K-approx}
   \widehat V_{n,\le K}
   =n^a\zeta_{n,a}Q_{n,K}^{(a)}+o_{L^1}(b_n).
\end{equation}

Finally, we use the common decomposition
\[
\begin{aligned}
\E\Big|
      \frac{\widehat V_n-n^aK_{a,\theta}\zeta_{n,a}}{b_n}
   \Big|&\le
   \frac{\E|\widehat V_{n,>K}|}{b_n}
   +\E\Big|
      \frac{\widehat V_{n,\le K}
            -n^a\zeta_{n,a}Q_{n,K}^{(a)}}{b_n}
   \Big|\\
   &\quad+
   \Big\|\frac{n^a\zeta_{n,a}}{b_n}\Big\|_2
   \left\|Q_{n,K}^{(a)}-K_{a,\theta}^{(K)}\right\|_2
   +|K_{a,\theta}^{(K)}-K_{a,\theta}|\,
    \frac{\E|n^a\zeta_{n,a}|}{b_n}.
\end{aligned}
\]
For fixed $K$, the second and third terms tend to zero as
$n\to\infty$ by \eqref{eq:V-K-approx} and
\eqref{eq:random-coeff-limit}, respectively, since
$n^a\zeta_{n,a}/b_n$ is bounded in $L^2$. By \eqref{eq:V-tail}, the
first term is bounded uniformly in $n$ by $C\sum_{k>K}\rho^{(k-1)/2}$, while
the fourth term is bounded by $C|K_{a,\theta}^{(K)}-K_{a,\theta}|$, because
$n^a\zeta_{n,a}/b_n$ is also bounded in $L^1$. Hence, taking
$\limsup_{n\to\infty}$ and then letting $K\to\infty$, we have
\[
   \widehat V_n=n^aK_{a,\theta}\zeta_{n,a}+o_{L^1}(b_n).
\]
Substituting the two choices of $\zeta_{n,a}$, we prove
\eqref{eq:V-critical} and \eqref{eq:V-strong}.
\end{proof}

\begin{proof}[Proof of Theorem~\ref{thm:fixed-main}]
The law of large numbers is Proposition~\ref{prop:fixed-lln}.
For the fluctuation limits, Lemma~\ref{lem:large-prime-reduction}
and the centered form of \eqref{eq:localization-transfer} allow us
to replace $T_n$ by $\widehat S_n$. Since $A_n^{(0)}$ is
deterministic and $\E\widehat U_n=0$,
by \eqref{eq:S-exact-decomposition} we have
\[
   \widehat S_n-\E\widehat S_n
   =\widehat U_n+\widehat V_n-\E\widehat V_n.
\]
Applying Proposition~\ref{prop:drift-projection}, using
$\E D_n=\E W_{a,\theta}=0$, and transferring the resulting expansion
back to $T_n$, we obtain
\[
   T_n-\E T_n=
   \begin{cases}
      \widehat U_n+o_P(b_n),&a<1/2,\\[2pt]
      \widehat U_n+K_{1/2,\theta}D_n+o_P(b_n),&a=1/2,\\[2pt]
      \widehat U_n+n^aK_{a,\theta}W_{a,\theta}+o_P(b_n),&a>1/2.
   \end{cases}
\]
The weak-memory limit follows from \eqref{eq:U-clt}. At
criticality, the two limiting Gaussian variables are independent by
\eqref{eq:joint-U-D}, so their variances add. For $a>1/2$, $\widehat U_n/n^a\toP0$ by
Proposition~\ref{prop:innovation-limit}, and the stated limit follows.
\end{proof}


\section{Proofs for variable retention}\label{sec:variable-proof}


In both boundary regimes, the genealogy is independent of the root labels,
but the label probability varies with $n$. We first analyze rare positive
families and then rare deletion families.


\subsection{Sparse retention}


Assume first that $\theta_n\sim\lambda/n$. We use the copy construction as
a random forest. Time $1$ is a root. For $j\ge2$, time $j$ starts a new
component with probability $1-a$ and otherwise attaches to a uniformly chosen
earlier vertex. Each root receives an independent Bernoulli$(\theta_n)$
label. A vertex is retained precisely when its root label is one.

Let $\mathcal R_n$ be the set of positive roots other than the initial root.
The indicators $\1_{\{j\in\mathcal R_n\}}$, $j\ge2$, are independent
Bernoulli variables with parameter $(1-a)\theta_n$. The point process $\sum_{j\in\mathcal R_n}\delta_{j/n}$
converges to a Poisson point process on $(0,1]$ with intensity
$(1-a)\lambda\,\dd u$, where $\delta_x$ denotes the unit point mass at $x$.
The Poisson convergence theorem for null arrays \cite[Theorem~30.1]{Kallenberg} applies because, on every interval
$(s,t]\subset(0,1]$, the sum of the success probabilities converges to
$(1-a)\lambda(t-s)$ and their maximum tends to zero. Each summand carries
at most one point, so the multiple-point condition in that theorem is
identically zero. The initial root is positive with probability $\theta_n\to0$ and has no effect on the limit.

For each $r\ge1$, let $\mathcal C_r(t)$ be the set of descendants of $r$,
including $r$, present at time $t$ if $r\le t$ and $r$ is a root;
otherwise, set $\mathcal C_r(t)=\varnothing$. Let
$C_r(t):=|\mathcal C_r(t)|$. We condition on $r\le t$ and $r$ being
a root. When $C_r(t)=c$, the next vertex joins this component with probability $ac/t$. Writing $c^{\overline m}=c(c+1)\cdots(c+m-1)$ for $m\ge1$ and $c^{\overline0}=1$, we obtain
\begin{equation}\label{eq:rising-factorial-recursion}
   \E\bigl(C_r(t+1)^{\overline m}\mid C_r(t)\bigr)
   =\left(1+\frac{am}{t}\right)C_r(t)^{\overline m}.
\end{equation}

\begin{proposition}\label{prop:Yule-marks}
If $r_n\in[n]$ and $r_n/n\to u\in(0,1]$, then, conditional on $r_n$ being a root,
\[
   C_{r_n}(n)\tod Z_u,
\]
where $Z_u$ has the geometric law \eqref{eq:Yule-mark}.
\end{proposition}

\begin{proof}
This is the standard discrete Yule family calculation; see also
\cite[Section~5.1]{BaurBertoinYule}. Iterating
\eqref{eq:rising-factorial-recursion}, we have
\begin{align}\label{eq:Yule-exponent-moment}
   \E\left(C_{r_n}(n)^{\overline m}
      \,\middle|\,r_n\text{ is a root}\right)
   =m!\prod_{t=r_n}^{n-1}\left(1+\frac{am}{t}\right)
   \longrightarrow m!u^{-am}.
\end{align}
These are the rising-factorial moments of a geometric variable with success
probability $u^a$. Tightness follows from the
case $m=1$. For each fixed $m$, applying \eqref{eq:Yule-exponent-moment} with $m+1$ in place of $m$
shows that $C_{r_n}(n)^{\overline m}$ is uniformly integrable. Thus
every subsequential limit in distribution has these moments. The
geometric law has a finite moment-generating function near zero and
is therefore moment determinate. Consequently, every subsequential
limit is $Z_u$, which proves the proposition.
\end{proof}

A marked Poisson limit also requires joint control over several families. The
next lemma establishes the asymptotic independence of finitely many families born
at times proportional to $n$.

\begin{lemma}\label{lem:Yule-joint}
Fix $s\ge1$. If distinct times $r_{n,1},\ldots,r_{n,s}\in[n]$
satisfy $r_{n,i}/n\to u_i\in(0,1]$, then, conditional on
$r_{n,1},\ldots,r_{n,s}$ all being roots,
\[
   \bigl(C_{r_{n,1}}(n),\ldots,C_{r_{n,s}}(n)\bigr)
   \tod (Z_{u_1},\ldots,Z_{u_s}),
\]
where $Z_{u_1},\ldots,Z_{u_s}$ are independent and each $Z_{u_i}$
has the geometric law \eqref{eq:Yule-mark}. The limiting locations
$u_i$ need not be distinct.
\end{lemma}

\begin{proof}
Fix nonnegative integers $m_1,\ldots,m_s$ and consider the
corresponding mixed rising-factorial moment. For fixed $n$, write
$r_i=r_{n,i}$, $\mathcal J=\{r_1,\ldots,r_s\}$, and
$r_*=\min_i r_i$. Let $I(t)=\{i:r_i\le t\}$ index the clusters that
have appeared by time $t$.

We first derive the one-step recursion. Away from the prescribed
root times, at most one active family grows. If its current size is
$c$, the identity $c\bigl((c+1)^{\overline m}-c^{\overline m}\bigr)=m c^{\overline m}$
shows that its growth increases the mixed rising-factorial moment by
$m/c$ times its current value. The new vertex joins this family with
probability $ac/t$. Summing over $i\in I(t)$, we obtain the one-step multiplier
$1+a\sum_{i\in I(t)}m_i/t$ for the conditional expectation given the genealogy through time $t$.
At a time when a prescribed new root is created, the existing
clusters do not grow, so the corresponding multiplier is one. The
birth of the $i$th cluster contributes the initial factor $m_i!$,
since its initial size is one. Root decisions at all other times are
independent of the conditioning event. Iterating the recursion from time $r_*$ to time $n$,
we obtain the exact conditional moment
\begin{equation}\label{eq:Yule-mixed-exact}
\begin{aligned}
 &\E\Big[\prod_{i=1}^s C_{r_i}(n)^{\overline{m_i}}
       \,\big|\,r_1,\ldots,r_s\text{ are roots}\Big]=\Big(\prod_{i=1}^s m_i!\Big)
   \prod_{\substack{r_*\le t<n\\t+1\notin\mathcal J}}
       \Big(1+\frac{a}{t}\sum_{i\in I(t)}m_i\Big).
\end{aligned}
\end{equation}

By the one-cluster moment formula \eqref{eq:Yule-exponent-moment}, the product of the individual
conditional moments is
\begin{align}\label{eq:Yule-mixed-individual}
   \prod_{i=1}^s
   \E\!\left(C_{r_i}(n)^{\overline{m_i}}
      \,\middle|\,r_i\text{ is a root}\right)
   =\Big(\prod_{i=1}^s m_i!\Big)
   \prod_{i=1}^s\prod_{t=r_i}^{n-1}
   \left(1+\frac{am_i}{t}\right).
\end{align}
We now compare the time product in \eqref{eq:Yule-mixed-exact}
with the double product in \eqref{eq:Yule-mixed-individual}.
At an ordinary time $t$, the product of the individual multipliers
active at that time satisfies
\[
   \prod_{i\in I(t)}\left(1+\frac{am_i}{t}\right)
   =1+\frac{a}{t}\sum_{i\in I(t)}m_i+O(t^{-2}).
\]
Uniformly when $r_*\ge\delta n$, all relevant times are at least $\delta n$, and the accumulated error from the ordinary times is
\[
   \sum_{t=\lceil\delta n\rceil}^{n-1}O(t^{-2})
   =O_{\delta,a,s,\boldsymbol m}(n^{-1}),
\]
where $\boldsymbol m=(m_1,\ldots,m_s)$.

At each prescribed root time $r_i$, we remove from the double product
in \eqref{eq:Yule-mixed-individual} the factors at $t=r_i-1$, one for
each cluster already active, since they have no counterpart in the
joint time product in \eqref{eq:Yule-mixed-exact}. Each such factor is $1+O((r_i-1)^{-1})$. Because
there are only finitely many prescribed root times, these removals
change the double product by a factor of $1+O(n^{-1})$.

Combining the two kinds of error, we find that, uniformly when
$r_*\ge\delta n$, the absolute logarithm of the ratio of the time
product in \eqref{eq:Yule-mixed-exact} to the double product in \eqref{eq:Yule-mixed-individual}
is $O_{\delta,a,s,\boldsymbol m}(n^{-1})$.
Consequently, for fixed $\delta>0$,
\begin{equation}\label{eq:Yule-mixed-uniform}
\begin{aligned}
 &\E\Big[\prod_{i=1}^s C_{r_i}(n)^{\overline{m_i}}
       \,\big|\,r_1,\ldots,r_s\text{ are roots}\Big]=\Big(\prod_{i=1}^s m_i!(n/r_i)^{am_i}\Big)
       \left(1+O_{\delta,a,s,\boldsymbol m}(n^{-1})\right),
\end{aligned}
\end{equation}
uniformly over distinct integers $r_i\in[\delta n,n]$. Here we use the product
estimate in Proposition~\ref{app:gamma-estimates} with $b=am_i$.
Thus every mixed rising-factorial moment factorizes asymptotically as the
product of the corresponding one-dimensional moments in
Proposition~\ref{prop:Yule-marks}. In particular, no separation between the limiting
locations $u_i$ is required.

Tightness follows from the marginal first-moment bounds in
\eqref{eq:Yule-mixed-uniform}. Moreover, applying
\eqref{eq:Yule-mixed-uniform} with $2m_i$ in place of $m_i$ bounds
the second moments of the corresponding mixed rising-factorial
products. Hence all products of fixed powers of the cluster sizes are
uniformly integrable. Since the limiting product law has a finite
moment-generating function near the origin, it is moment
determinate and the joint method of moments now completes the proof.
\end{proof}

The marked point process argument applies away from time zero. We
next bound the expected contribution of families with very early
positive roots, so that they can be removed when the truncation
parameter tends to zero.

\begin{lemma}\label{lem:early-positive-roots}
For $0<\delta<1$, let $M_n^{\le\delta}$ be the number of retained
vertices whose positive root time is at most $\delta n$. Then
\begin{equation}\label{eq:early-root-bound}
   \limsup_{n\to\infty}\E M_n^{\le\delta}
   \le C\delta^{1-a},
\end{equation}
where $C$ is independent of $n$ and $\delta$.
\end{lemma}

\begin{proof}
We treat the initial root separately. Its contribution is at most
\[
\theta_n\E C_1(n)=O(\theta_n n^a)=o(1).
\]
For $r\ge2$, the probability that $r$ is a
positive root is $(1-a)\theta_n$, and the root label is independent of the forest. By the moment recursion \eqref{eq:rising-factorial-recursion} and Proposition \ref{app:gamma-estimates}, we have
\[
   \E\left(C_r(n)\mid r\text{ is a root}\right)
   =\prod_{j=r}^{n-1}\Big(1+\frac aj\Big)
   \le C(n/r)^a.
\]
Summing over the roots at most $\delta n$, we obtain
\[
   \E M_n^{\le\delta}
   \le o(1)+C\theta_n\sum_{r\le\delta n}(n/r)^a \le C\delta^{1-a}+o(1),
\]
which completes the proof.
\end{proof}

\begin{proof}[Proof of \eqref{eq:sparse-count}]
We first prove a marked Poisson limit for roots born after time
$\delta n$. We then sum their marks by truncating large marks, and
finally remove the time truncation by controlling the early roots in
both the prelimit and the limit.

Fix $\delta\in(0,1)$ and introduce the marked point process
\[
   \Xi_n^{\delta}:=
   \sum_{\substack{r\in\mathcal R_n\\r>\delta n}}
   \delta_{(r/n,C_r(n))}
   \quad\text{on }[\delta,1]\times\N.
\]
We equip the space of locally finite point measures on
$[\delta,1]\times\N$ with the vague topology, where $\N$ has the discrete
topology. Define the intensity measure
\[
   \nu_\delta(\dd u,\dd z)
   :=(1-a)\lambda\,\dd u\,K_u(\dd z),
   \qquad
   K_u(\{z\}):=\Pp(Z_u=z),\quad z\in\N.
\]

We first compute the factorial moments of $\Xi_n^\delta$. Fix
$m\ge1$ and let $f_1,\ldots,f_m\in C_c^+([\delta,1]\times\N)$. Extend each $f_i$ to mark zero by setting $f_i(u,0)=0$. Independence of the root decisions and their labels yields
\begin{equation}\label{eq:sparse-factorial-prelimit}
\begin{aligned}
 &\E\sum_{\substack{\max\{1,\delta n\}<r_1,\ldots,r_m\le n\\r_1,\ldots,r_m\text{ distinct}}}
   \prod_{i=1}^m f_i(r_i/n,C_{r_i}(n))
   \prod_{i=1}^m\1_{\{r_i\in\mathcal R_n\}}\\
 &\qquad=((1-a)\theta_n)^m
   \sum_{\substack{\max\{1,\delta n\}<r_1,\ldots,r_m\le n\\r_1,\ldots,r_m\text{ distinct}}}
   \E\Big[
      \prod_{i=1}^m f_i(r_i/n,C_{r_i}(n))
      \,\big|\, r_1,\ldots,r_m\text{ are roots}
   \Big].
\end{aligned}
\end{equation}
The conditioning event only requires the displayed times to be
roots. Their positive labels are independent of the forest and do not
change the conditional law of the cluster sizes.

By Lemma~\ref{lem:Yule-joint}, if $r_i/n\to u_i$ for each $i$, then
the conditional expectations in
\eqref{eq:sparse-factorial-prelimit} converge to $\prod_{i=1}^m\E f_i(u_i,Z_{u_i})$. Moreover, the convergence is uniform over distinct root tuples with
$r_i/n\in[\delta,1]$, in the sense that
\[
   \E\!\Big[
      \prod_{i=1}^m f_i(r_i/n,C_{r_i}(n))
      \,\big|\,r_1,\ldots,r_m\text{ are roots}
   \Big]
   -\prod_{i=1}^m\E f_i(r_i/n,Z_{r_i/n})
   \longrightarrow0
\]
uniformly over such tuples. Indeed, suppose it
fails. After taking a subsequence, we may assume that
$r_i/n\to u_i\in[\delta,1]$ for every $i$, since $[\delta,1]^m$ is
compact. Lemma \ref{lem:Yule-joint} and the continuity of
$u\mapsto\E f_i(u,Z_u)$ then yield a contradiction. The continuity
of these functions follows from the finite mark support of $f_i$.
Since $(1-a)\theta_n\sim(1-a)\lambda/n$, a Riemann sum argument
therefore identifies the limit of the right-hand side of
\eqref{eq:sparse-factorial-prelimit} as
\[
   \prod_{i=1}^m
   (1-a)\lambda\int_\delta^1\E f_i(u,Z_u)\,\dd u
   =\prod_{i=1}^m\int_{[\delta,1]\times\N}f_i\,\dd\nu_\delta.
\]
The Riemann sum also contains the diagonal tuples omitted from the
factorial sum, but their total contribution is $O_m(n^{-1})$,
because there are $O_m(n^{m-1})$ such tuples and
$\theta_n^m=O(n^{-m})$. Thus all factorial moments of
$\Xi_n^\delta$ converge to those of a Poisson process with intensity
$\nu_\delta$.

Let
\[
   B_n^\delta:=\#\{r\in\mathcal R_n:r>\delta n\},
\]
and then
\[
   \sup_n\E 3^{B_n^\delta}
   \le\sup_n\bigl(1+2(1-a)\theta_n\bigr)^n<\infty.
\]
Since $\Xi_n^\delta(B)\le B_n^\delta$ for every
$B\subset [\delta,1]\times\N$, Proposition~\ref{app:poisson-factorial-criterion}
now yields
\begin{equation}\label{eq:sparse-marked-process-limit}
   \Xi_n^\delta\tod\Xi^\delta:=
   \sum_{u\in\mathcal R\cap(\delta,1]}\delta_{(u,Z_u)}
\end{equation}
in the vague topology.

We next pass from the point processes to the sums of their marks.
For $K\ge1$, define
\[
   F_K(\xi):=\int z\1_{\{z\le K\}}\,\xi(\dd u,\dd z).
\]
The integrand is continuous with compact support on
$[\delta,1]\times\N$, so $F_K$ is continuous for the vague topology. Hence,
$F_K(\Xi_n^\delta)\tod F_K(\Xi^\delta)$, for each fixed $K$.

It remains to control the marks above $K$. By the moment recursion \eqref{eq:rising-factorial-recursion}
and the product estimates \eqref{eq:app-gamma-product}, we obtain, uniformly for $r>\delta n$,
\[
   \E\left(C_r(n)^2\mid r\text{ is a root}\right)
   \le C_\delta.
\]
Consequently, with $\alpha_n=(1-a)\theta_n$, we have
\[
\begin{aligned}
 &\alpha_n\sum_{\delta n<r\le n}
 \E\left(C_r(n)\1_{\{C_r(n)>K\}}
       \,\middle|\,r\text{ is a root}\right)
 \le\frac{C_\delta}{K},
\end{aligned}
\]
where we used that $n\theta_n$ is bounded. This bounds the expected
mark-truncation error uniformly in $n$. For the limiting process,
the marked Campbell formula \cite[Proposition~2.7 and Theorem~5.6]{LastPenrose} yields
\[
   \E\sum_{u\in\mathcal R\cap(\delta,1]}
      Z_u\1_{\{Z_u>K\}}=\lambda(1-a)\int_\delta^1\E(Z_u\1_{\{Z_u>K\}})\dd u
   \le\frac{C_\delta}{K},
\]
since $\sup_{u\ge\delta}\E Z_u^2<\infty$. We can therefore apply the converging-together theorem with $K\to\infty$ and
conclude that
\begin{equation}\label{eq:sparse-truncated-count-limit}
   \sum_{\substack{r\in\mathcal R_n\\r>\delta n}}C_r(n)
   \tod
   N_{a,\lambda}^{>\delta}
   :=\sum_{u\in\mathcal R\cap(\delta,1]}Z_u.
\end{equation}

Finally, we remove the time truncation. For $n\ge\delta^{-1}$, the
difference between $M_n$ and the left-hand side of
\eqref{eq:sparse-truncated-count-limit} is precisely
$M_n^{\le\delta}$. By Lemma~\ref{lem:early-positive-roots} and
Markov's inequality, this difference is negligible when we first let
$n\to\infty$ and then $\delta\downarrow0$. The corresponding part of
the limiting process is also negligible, since, by Campbell's formula
again,
\[
   \E\sum_{u\in\mathcal R\cap(0,\delta]}Z_u
   =(1-a)\lambda\int_0^\delta u^{-a}\,\dd u
   =\lambda\delta^{1-a}.
\]

This completes the proof of \eqref{eq:sparse-count}.
\end{proof}

We now transfer the limit for the retained count to the logarithmic least
common multiple.

\begin{lemma}\label{lem:sparse-log-sum}
In the sparse regime,
\begin{equation}\label{eq:sparse-log-sum}
   \frac1{\log n}\sum_{1\le j\le n}X_{n,j}\log j-M_n\toP0.
\end{equation}
\end{lemma}

\begin{proof}
By the definition of $M_n$ and $\E X_{n,j}=\theta_n$, we have
\[
\begin{aligned}
 \E\Big|\frac1{\log n}\sum_{1\le j\le n}X_{n,j}\log j-M_n\Big|
   &=\frac{\theta_n}{\log n}\sum_{j=1}^n\log(n/j)\le\frac{\theta_n(n-1)}{\log n}\longrightarrow0,
\end{aligned}
\]
where the last inequality follows from $\log(n!)\ge n\log n-n+1$. Thus the convergence holds in $L^1$,
and hence in probability.
\end{proof}

\begin{lemma}\label{lem:sparse-gcd}
In the sparse regime,
\begin{equation}\label{eq:sparse-gcd}
   \frac1{\log n}\sum_{1\le i<j\le n}
   X_{n,i}X_{n,j}\log\gcd(i,j)\toP0.
\end{equation}
\end{lemma}

\begin{proof}
Fix $K$ and $0<\eta<1$. We first prove that
\begin{align}\label{eq:sparse-gcd-eta}
   \Pp\bigl(M_n\le K,\ \exists\,i<j
   \text{ retained with }\gcd(i,j)>n^\eta\bigr)
   \longrightarrow0.
\end{align}
Let us begin with pairs from different roots. Conditional on the forest,
their root labels are independent, so each specified pair is retained
with probability at most $\theta_n^2$. Summing over all such pairs
with $\gcd(i,j)>n^\eta$, we obtain the bound
$\theta_n^2\sum_{h>n^\eta}\left\lfloor n/h\right\rfloor^2\le Cn^{-\eta}$.

We next consider pairs from the same positive cluster. Let
$A_r=\{r\text{ is a positive root}\}$, and let $\mathcal H_t^*$ be
the $\sigma$-field generated by the forest and root labels through
time $t$. On $A_r\cap\{M_n\le K\}$, we have $C_r(t)\le K$ for every
$t\le n$. For $r<i<j$, we enlarge the event by replacing $M_n\le K$ with
$C_r(j-1)\le K$. Conditioning successively on
$\mathcal H_{j-1}^*$ and $\mathcal H_{i-1}^*$, we therefore obtain
\[
\begin{aligned}
 &\Pp\bigl(A_r,\ i,j\in\mathcal C_r(n),\ M_n\le K\bigr)\\
 &\le\E\Big[\1_{A_r}\1_{\{i\in\mathcal C_r(i)\}}
          \1_{\{C_r(j-1)\le K\}}\frac{aC_r(j-1)}{j-1}\Big]\\
 &\le \frac{aK}{j-1}
   \Pp\bigl(A_r,\ i\in\mathcal C_r(i),\ C_r(i-1)\le K\bigr)\le \frac{(aK)^2}{(i-1)(j-1)}\Pp(A_r).
\end{aligned}
\]
For the roots $2\le r<i$, we have
$\Pp(A_r)=(1-a)\theta_n$, so their total contribution is at most
\[
   (1-a)\theta_n(i-2)\frac{(aK)^2}{(i-1)(j-1)}
   \le \frac{C_K\theta_n}{j}.
\]
If $r=i$, then $i$ itself is a positive root, and only $j$ needs to
join the cluster. This case contributes at most
$(1-a)\theta_n aK/(j-1)$. For $r=1$, the initial root is positive
with probability $\theta_n$, and both $i$ and $j$ must join its
cluster. Its contribution is at most
$C_K\theta_n/\{(i-1)(j-1)\}$. Combining these cases, we find that
\begin{equation}\label{eq:same-root-pair}
   \Pp(i,j\text{ are retained from the same root},\ M_n\le K)
   \le \frac{C_K\theta_n}{j}.
\end{equation}
The case $i=1$ contributes zero because $\gcd(1,j)=1$.
It remains to sum over $i\ge2$ with $\gcd(i,j)>n^\eta$. Writing
$i=uh$ and $j=vh$ with $u<v$, we have
\[
   \sum_{\substack{2\le i<j\le n\\ \gcd(i,j)>n^\eta}}\frac1j
   \le \sum_{h>n^\eta}\sum_{v\le n/h}\frac{v-1}{vh}
   \le C\sum_{h>n^\eta}\frac{n}{h^2}\le Cn^{1-\eta}.
\]
Together with \eqref{eq:same-root-pair}, this yields
\[
   \sum_{\substack{1\le i<j\le n\\ \gcd(i,j)>n^\eta}}
   \Pp(i,j\text{ are retained from the same root},\ M_n\le K)
   \le C_Kn^{-\eta}.
\]
Combining this with the different-root estimate, we prove \eqref{eq:sparse-gcd-eta}.

Note that, on $\{M_n\le K\}$, there are at most $K^2$ retained pairs. Outside
the exceptional event, each of them satisfies
$\log\gcd(i,j)\le\eta\log n$, so the normalized sum in
\eqref{eq:sparse-gcd} is at most $K^2\eta$.
Given $\varepsilon,\delta>0$, we first choose $K$ large enough that
\[
   \sup_n\Pp(M_n>K)
   \le\frac{\sup_n n\theta_n}{K}<\delta,
\]
using Markov's inequality and $\E M_n=n\theta_n$. We then choose
$0<\eta<1$ small enough that $K^2\eta<\varepsilon$. Since the
exceptional event has probability tending to zero, the limsup of the
probability that the left-hand side of \eqref{eq:sparse-gcd} exceeds
$\varepsilon$ is at most $\delta$. Letting $\delta\downarrow0$, we prove the lemma.
\end{proof}

\begin{proof}[Proof of Theorem~\ref{thm:sparse-main}]
We first prove a deterministic comparison between a product and the
least common multiple. The empty and singleton cases are immediate.
For $m\ge2$ and positive integers $j_1,\ldots,j_m$, we claim that
\begin{align}\label{eq:sparse-product-exponent}
   \prod_{r=1}^m j_r\Big/\prod_{r<s}\gcd(j_r,j_s)
   \le \LCM(j_1,\ldots,j_m)
   \le\prod_{r=1}^m j_r.
\end{align}
The upper bound is immediate. For the lower bound, fix a prime and
order its exponents in $j_1,\ldots,j_m$ as
$e_1\le\cdots\le e_m$. Its exponent on the left-hand side is
\[
   \sum_{r=1}^m e_r-\sum_{r<s}\min(e_r,e_s)
   =e_m-\sum_{r=1}^{m-2}(m-r-1)e_r
   \le e_m,
\]
where $e_m$ is its exponent in the least common multiple.
This proves \eqref{eq:sparse-product-exponent} prime by prime.

Applying this inequality to the retained integers and taking
logarithms, we obtain
\[
   0\le \sum_{1\le j\le n}X_{n,j}\log j-T_n
   \le\sum_{1\le i<j\le n}
      X_{n,i}X_{n,j}\log\gcd(i,j).
\]
Lemmas~\ref{lem:sparse-log-sum} and \ref{lem:sparse-gcd} therefore imply
\[
   \frac{T_n}{\log n}-M_n\toP0.
\]
Combining this with \eqref{eq:sparse-count}, we obtain \eqref{eq:sparse-lcm}.
The proof of Theorem~\ref{thm:sparse-main} is complete.
\end{proof}


\subsection{Near-complete retention}


Assume now that
\[
   q_n:=1-\theta_n\sim\frac{\lambda\log n}{n}.
\]
Set $Y_{n,j}=1-X_{n,j}$. The deleted row obeys the same model with
root label probability $q_n$. We call the primes in $\calP_n^+=\calP\cap(n/2,n]$ upper primes.
Let $\mathcal R_n^0$ be the set of genealogical roots, including the
initial root. For $r\in\mathcal R_n^0$, define $N_{n,r}=\#\bigl(\mathcal C_r(n)\cap\calP_n^+\bigr)$,
and set $N_{n,r}=0$ when $r\notin\mathcal R_n^0$. Conditional on the
forest, write $B_{n,r}$ for the independent Bernoulli$(q_n)$ deletion label
of root $r$. Each upper prime belongs to exactly one root cluster
and is deleted precisely when its root label is one. Hence the number
of deleted upper primes has the exact representation
\begin{equation}\label{eq:deleted-root-representation}
   Z_n:=\sum_{p\in\calP_n^+}Y_{n,p}
      =\sum_{r\in\mathcal R_n^0}B_{n,r}N_{n,r}.
\end{equation}

In the next lemma, we estimate the first moment and the second
factorial moment of $N_{n,r}$, conditional on $r$ being a root. The
convention $N_{n,r}=0$ for nonroots allows us to use these estimates
in unconditional sums by multiplying by
$\Pp(r\in\mathcal R_n^0)$, which is $1$ for $r=1$ and $1-a$ for
$r\ge2$.

\begin{lemma}\label{lem:upper-prime-cluster}
Conditional on $r$ being a root, there is a constant $C<\infty$ such
that the following bounds hold: for $1\le r\le n/2$,
\begin{align}
   \E N_{n,r}
   &\le C\frac{n^a}{r^a\log n},
   \label{eq:upper-first}\\
   \E[N_{n,r}(N_{n,r}-1)]
   &\le C\frac{n^{2a}}{r^{2a}\log^2 n};
   \label{eq:upper-second}
\end{align}
for $n/2<r\le n$,
\begin{equation}\label{eq:upper-late-second}
   \E[N_{n,r}(N_{n,r}-1)]
   \le C\Big(\frac{\1_{\{r\in\calP\}}}{\log n}
   +\frac1{\log^2 n}\Big).
\end{equation}
\end{lemma}

\begin{proof}
All expectations in this proof are conditional on $r$ being a root.
By using the rising-factorial recursion
\eqref{eq:rising-factorial-recursion} and the gamma-ratio bound \eqref{eq:app-gamma-product}, we have
\[
   \E C_r(t)\le C(t/r)^a
   \quad\text{and}\quad
   \E C_r(t)^2\le C(t/r)^{2a}.
\]
If $r<p$, conditioning on the genealogy through time $p-1$ shows
that $p$ joins the cluster with probability
\[
   \frac{a}{p-1}\E C_r(p-1)
   \le Cr^{-a}p^{a-1}.
\]
Summing over the upper primes, we get
\[
   \E N_{n,r}
   \le Cr^{-a}\sum_{p\in\calP_n^+}p^{a-1}
   \le \frac{C n^a}{r^a\log n},
\]
where the last step uses partial summation and the prime number
theorem. This proves \eqref{eq:upper-first}.

Let $\mathcal H_s$ be the $\sigma$-field generated by the genealogy through
time $s$. For $r<p_1<p_2$, the conditional mean evolution of the cluster is
\[
   \E\bigl(C_r(t)\mid\mathcal H_s\bigr)
   =C_r(s)\prod_{h=s}^{t-1}\left(1+\frac ah\right),
   \qquad r\le s<t.
\]
Conditioning successively at $p_2-1$ and $p_1$, we obtain
\[
\begin{aligned}
 &\Pp(p_1,p_2\text{ join the cluster})=\frac{a}{p_2-1}
   \prod_{h=p_1}^{p_2-2}\Big(1+\frac ah\Big)
   \E\left[\1_{\{p_1\text{ joins}\}}C_r(p_1)\right].
\end{aligned}
\]
The last expectation is $\frac{a}{p_1-1}\E C_r(p_1-1)^{\overline 2}$.
Applying the rising-factorial recursion and the gamma-ratio bound again, we obtain
\[
   \E C_r^{\overline 2}(p_1-1)\le C(p_1/r)^{2a}
   \quad\text{and}\quad
   \prod_{h=p_1}^{p_2-2}\Big(1+\frac ah\Big)
      \le C(p_2/p_1)^a.
\]
It follows that the probability that $p_1$ and $p_2$ both join the
cluster is at most $Cr^{-2a}p_1^{a-1}p_2^{a-1}$. Since
\[
   N_{n,r}(N_{n,r}-1)
   =2\sum_{\substack{p_1,p_2\in\calP_n^+\\p_1<p_2}}
      \1_{\{p_1,p_2\text{ join the cluster}\}},
\]
we conclude that
\[
   \sum_{\substack{p_1,p_2\in\calP_n^+\\p_1<p_2}}
   \Pp(p_1,p_2\text{ join the cluster})
   \le Cr^{-2a}
   \bigg(\sum_{p\in\calP_n^+}p^{a-1}\bigg)^2,
\]
which proves \eqref{eq:upper-second}.

For $n/2<r\le n$, let
\[
   D_{n,r}:=\#\{p\in\calP_n^+:p>r,\ p\text{ belongs to the cluster of }r\}.
\]
In this range, $p/r\le2$, so the attachment bounds above become
$O(1/p)$ for one vertex and $O(1/(p_1p_2))$ for two vertices. Summing over primes and applying
the prime number theorem, we get
\begin{align}\label{eq:sparse-complete-1}
   \E D_{n,r}\le C\sum_{\substack{ p\in\calP\\r<p\le n}}\frac1p
      =O\left(\frac1{\log n}\right) \quad\text{and}\quad
      \E[D_{n,r}(D_{n,r}-1)]
      =O\left(\frac1{\log^2 n}\right),
\end{align}
uniformly in $n/2<r\le n$.

To relate these bounds to $N_{n,r}$, note that
\[
   N_{n,r}=\1_{\{r\in\calP\}}+D_{n,r},
\]
and hence
\[
   N_{n,r}(N_{n,r}-1)
   =D_{n,r}(D_{n,r}-1)+2\1_{\{r\in\calP\}}D_{n,r}.
\]
Taking expectations and using \eqref{eq:sparse-complete-1}, we obtain
\eqref{eq:upper-late-second}.
\end{proof}

These bounds allow us to discard early roots and to ignore, with high
probability, deleted families containing more than one upper prime.
The remaining count has a Poisson limit.

\begin{lemma}\label{lem:deleted-primes-poisson}
The number $Z_n$ of deleted upper primes, defined in
\eqref{eq:deleted-root-representation}, satisfies $Z_n\tod\Pi(\lambda/2)$.
\end{lemma}

\begin{proof}
Let $\ell_n=n/\log n$. We need two truncation estimates: the contribution from early roots
must be negligible, and late clusters must rarely contain more than
one deleted upper prime. In the sums below, expectations are
unconditional and $N_{n,r}=0$ when $r$ is not a root. We can thus use
the conditional bounds in Lemma~\ref{lem:upper-prime-cluster} after
multiplying by the root probability, which we bound by one. We claim
that
\begin{equation}\label{eq:early-upper-negligible}
   q_n\E\sum_{r\le\ell_n}N_{n,r}\longrightarrow0
\end{equation}
and
\begin{equation}\label{eq:upper-pair-negligible}
   q_n\E\sum_{\ell_n<r\le n}N_{n,r}(N_{n,r}-1)\longrightarrow0.
\end{equation}

For \eqref{eq:early-upper-negligible}, by the first-moment bound \eqref{eq:upper-first}, we have
\[
   q_n\E\sum_{r\le\ell_n}N_{n,r}
   \le \frac{Cq_n n^a}{\log n}\sum_{r\le\ell_n}r^{-a}.
\]
Since $q_n\sim\lambda(\log n)/n$ and
$\sum_{r\le\ell_n}r^{-a}=O(\ell_n^{1-a})$,  the right-hand side is
$O((\ell_n/n)^{1-a})$ and hence tends to zero.

For \eqref{eq:upper-pair-negligible}, by \eqref{eq:upper-second}, the roots
between $\ell_n$ and $n/2$ contribute at most
\[
   \frac{Cq_n n^{2a}}{\log^2n}
   \sum_{\ell_n<r\le n/2}r^{-2a}
   =
   \begin{cases}
      O(1/\log n),&a<1/2,\\[2pt]
      O(\log\log n/\log n),&a=1/2,\\[2pt]
      O((\log n)^{2a-2}),&a>1/2.
   \end{cases}
\]
Each bound tends to zero
because $a<1$. By \eqref{eq:upper-late-second}, the roots larger than
$n/2$ contribute at most
\[
   Cq_n\Big(\frac{\pi(n)-\pi(n/2)}{\log n}
                 +\frac{n}{\log^2n}\Big)
   =O(1/\log n).
\]
Combining the two parts, we obtain \eqref{eq:upper-pair-negligible}.

Let $\mathcal G_n$ be the $\sigma$-field generated by the genealogical
forest, and define
\[
   Z_n^>:=\sum_{\substack{r\in\mathcal R_n^0\\\ell_n <r\le n}}
      B_{n,r}N_{n,r}.
\]
By Markov's inequality and \eqref{eq:early-upper-negligible},
\begin{align}\label{eq:sparse-complete-approx}
   \Pp(Z_n\ne Z_n^>)
   \le q_n\E\sum_{r\le\ell_n}N_{n,r}\longrightarrow0.
\end{align}
We now compute the conditional probability generating function of
$Z_n^>$. Given $\mathcal G_n$, the variables $B_{n,r}$ are
independent, so for $0\le z\le1$,
\[
   \E(z^{Z_n^>}\mid\mathcal G_n)
   =\prod_{\substack{r\in\mathcal R_n^0\\ \ell_n <r\le n}}
      \left\{1+q_n(z^{N_{n,r}}-1)\right\}.
\]
For $q_n\le1/2$, we apply
$|\log(1+x)-x|\le2x^2$ with
$x=q_n(z^{N_{n,r}}-1)$. Since $|z^{N_{n,r}}-1|\le1$ and there are at
most $n$ factors, it follows that
\[
   \log\E(z^{Z_n^>}\mid\mathcal G_n)
   =q_n\sum_{\ell_n<r\le n}(z^{N_{n,r}}-1)+O(nq_n^2).
\]
The last error tends to zero because $nq_n^2=O(\log^2n/n)$.
Moreover,
\[
   |z^m-1-(z-1)m|
   \le\frac12(1-z)^2m(m-1),
\]
so the pair estimate \eqref{eq:upper-pair-negligible} implies
\[
   q_n\sum_{\ell_n<r\le n}(z^{N_{n,r}}-1)
   =(z-1)q_n\sum_{\ell_n<r\le n}N_{n,r}+o_P(1).
\]
Every prime in $\calP_n^+$ belongs to exactly one root cluster. Using
\eqref{eq:early-upper-negligible} and the prime number theorem, we
therefore obtain
\[
   q_n\sum_{\ell_n<r\le n}N_{n,r}
   =q_n|\calP_n^+|+o_P(1)
   \toP\frac\lambda2.
\]
Thus the conditional probability generating functions of $Z_n^>$ converge in probability to
$\exp\{(\lambda/2)(z-1)\}$. Taking expectations and using the bounded convergence theorem, we prove the Poisson limit for $Z_n^>$. The same limit holds for $Z_n$ by \eqref{eq:sparse-complete-approx}.
\end{proof}

\begin{remark}\label{rem:near-moments}
For $a\ge1/2$, the truncation by $\ell_n$ is needed in
\eqref{eq:upper-pair-negligible}. If we sum the bound
\eqref{eq:upper-second} over all roots instead, the resulting error
is only $O(1)$ at $a=1/2$ and $O(n^{2a-1}/\log n)$ for $a>1/2$,
neither of which tends to zero. Restricting the sum to
$r>\ell_n$ makes the error tend to zero. For $a<1/2$, no truncation
is needed for this estimate.
\end{remark}

\begin{proof}[Proof of Theorem~\ref{thm:near-complete}]
The prime-power decomposition follows the strategy of
\cite[Section~6.2]{AKM}. The new point is to control the probability
that the required deletions share a genealogical root.
By the von Mangoldt representation \eqref{eq:intro-vm}, we have
\[
   \Delta_n=\sum_{m\le n}\Lambda(m)(1-I_n(m)).
\]
For $2\le m\le n/2$, the event $I_n(m)=0$ implies that both $m$ and
$2m$ are deleted. The value $m=1$ is irrelevant since
$\Lambda(1)=0$. The genealogy formula in Lemma~\ref{lem:genealogy}, applied with root label
probability $q_n$, therefore implies
\[
   \Pp(I_n(m)=0)
   \le\Pp(Y_{n,m}=Y_{n,2m}=1)
   =q_n^2+q_n(1-q_n)\Pp(m\sim2m).
\]
The estimates in Proposition~\ref{app:arithmetic-estimates}, together with \eqref{eq:root-probability-bounds} and $q_n\sim\lambda(\log n)/n$, then yield that
\[
   \sum_{m\le n/2}\Lambda(m)\Pp(I_n(m)=0)\longrightarrow0.
\]
By Markov's inequality, the total contribution to $\Delta_n$ from
$m\le n/2$ thus converges to zero in probability.

For every $m\in(n/2,n]$, the integer $m$ is the only multiple of itself in
$[n]$, so $1-I_n(m)=Y_{n,m}$. Hence the expected contribution from prime
powers $p^r\in(n/2,n]$ with $r\ge2$ is at most
\[
   q_n\sum_{\substack{p^r\le n\\r\ge2}}\log p
   =q_n\bigl(\psi(n)-\vartheta(n)\bigr)=o(1),
\]
where we used $\psi(n)-\vartheta(n)=O(\sqrt n+n^{1/3}\log n)$. Markov's inequality
therefore makes this remainder negligible in probability.

After these reductions, the only remaining terms are those from the
upper primes, which satisfy the same identity $1-I_n(p)=Y_{n,p}$. We have thus shown that
\[
   \frac{\Delta_n}{\log n}
   -\frac1{\log n}\sum_{p\in\calP_n^+}\log p\,Y_{n,p}
   \toP0.
\]
Since $\log(n/2)\le\log p\le\log n$ on $\calP_n^+$,
\[
   \Big|
   \frac1{\log n}\sum_{p\in\calP_n^+}\log p\,Y_{n,p}-Z_n
   \Big|
   \le\frac{\log 2}{\log n}Z_n.
\]
Because $\E Z_n=q_n|\calP_n^+|\to\lambda/2$, this error tends to
zero in $L^1$. Lemma~\ref{lem:deleted-primes-poisson} and Slutsky's
theorem now complete the proof.
\end{proof}


\section{Concluding remarks}\label{sec-final}


Two natural directions remain open. The first concerns
weighted divisibility functionals
\[
   \mathcal W_n(w):=\sum_{m\le n}w(m)
      \1_{\{A_n\cap m\N\ne\varnothing\}},
\]
where $w$ is an arithmetic weight. We have treated $w=\Lambda$. For an
indeterminate $q$, define $[j]_q=(q^j-1)/(q-1)$. The degree of the polynomial
LCM of $[j]_q$, $j\in A_n$, equals $\mathcal W_n(w)$ with $w(1)=0$ and
$w(d)=\varphi(d)$ for $d\ge2$, where $\varphi$ is Euler's totient function.
Related questions for random $q$-integers under independent selection are
studied in \cite{Sanna}. Determining the fluctuation scale for other weights would clarify how memory interacts with the arithmetic features of the weight.

The second is a functional limit theorem at fixed density.
Such a result would describe the evolution of the process
$(T_{\lfloor nt\rfloor})_{0\le t\le1}$ with suitable centering and
normalization. Extending the proof would require a count projection that is uniform in time and suitable
tightness estimates. These are the main additional steps beyond the endpoint
limits proved here.


\appendix


\section{Auxiliary estimates and a Poisson criterion}\label{app:auxiliary-estimates}


We collect the arithmetic estimates used repeatedly in the proofs.  The first two propositions state standard analytic facts, while
the final proposition states the precise Poisson criterion needed for the
argument based on marked roots.

\begin{proposition}[Prime-number estimates]\label{app:arithmetic-estimates}
As $x\to\infty$,
\begin{equation}\label{eq:app-pnt}
   \vartheta(x)\sim x,
   \qquad
   \pi(x)\sim\frac{x}{\log x}.
\end{equation}
In particular, $\vartheta(x)=O(x)$ and $\psi(x)=O(x)$. Moreover,
\begin{gather}
   \sum_{\substack{p^r\le x\\r\ge2}}\log p
   =O(\sqrt x+x^{1/3}\log x),
   \label{eq:app-prime-powers}\\
   \sum_{m\le x}\frac{\Lambda(m)}m=O(\log x),
   \qquad
   \sum_{m\le x}\frac{\Lambda(m)\log m}{m}=O(\log^2x),
   \label{eq:app-mangoldt-harmonic}\\
   \sum_{m\le x}\Lambda(m)m^{\beta-1}
   =O_\beta(x^\beta),
   \qquad \beta>0.
   \label{eq:app-mangoldt-power}
\end{gather}
For each fixed $a\in[0,1)$,
\begin{equation}\label{eq:app-upper-primes}
   \sum_{x/2<p\le x}p^{a-1}=O_a\left(\frac{x^a}{\log x}\right),
   \qquad
   \sum_{x/2<p\le x}\frac1p=O\left(\frac1{\log x}\right).
\end{equation}
\end{proposition}

These estimates are standard consequences of the prime number theorem,
Chebyshev's bounds, and Abel summation; see
\cite[Chapters~4 and~13]{Apostol}.

\begin{proposition}[Gamma ratios and power sums]\label{app:gamma-estimates}
For fixed $a\in[0,1)$,
\[
   \gamma_n=\frac{\Gamma(n+a)}{\Gamma(1+a)\Gamma(n)}
   \sim\frac{n^a}{\Gamma(1+a)}.
\]
For every $\delta\in(0,1)$, both
\[
   \frac{\gamma_j}{\gamma_n}-\left(\frac jn\right)^a
   \quad\text{and}\quad
   \frac{\gamma_j}{j^a}-\frac1{\Gamma(1+a)}
\]
converge to zero uniformly for $\delta n\le j\le n$. Furthermore,
\begin{equation}\label{eq:app-gamma-sums}
   \sum_{j=1}^n\gamma_j^{-2}\sim
   \begin{cases}
      \displaystyle
      \frac{\Gamma(1+a)^2}{1-2a}\,n^{1-2a},&a<1/2,\\[2mm]
      \Gamma(3/2)^2\log n,&a=1/2,
   \end{cases}
\end{equation}
while the series converges when $a>1/2$.
For every fixed $b\ge0$ and integers $1\le r\le t$,
\begin{equation}\label{eq:app-gamma-product}
   \prod_{h=r}^{t-1}\Big(1+\frac bh\Big)
      =\frac{\Gamma(t+b)\Gamma(r)}{\Gamma(t)\Gamma(r+b)}
      \asymp_b (t/r)^b,
\end{equation}
where $\asymp_b$ means that the ratio is bounded above and below by
positive constants depending only on $b$. If $r\ge\delta t$ with
$\delta>0$ fixed, the last comparison can be sharpened to
$(t/r)^b(1+O_{\delta,b}(t^{-1}))$ as $t\to\infty$.
\end{proposition}

The gamma-ratio assertions follow from the uniform form of Stirling's
formula; see \cite[Sec.~5.11(iii)]{NIST}. The remaining statements follow from
the standard asymptotics of power sums.

\begin{proposition}[A factorial-moment Poisson criterion]
\label{app:poisson-factorial-criterion}
Let $E$ be a locally compact, second-countable Hausdorff space, and let $\nu$
be a Radon measure on $E$. Let $\Xi_n$ be point processes on $E$. Suppose
that, for every
$m\ge1$, $f_1,\ldots,f_m\in C_c^+(E)$, where $C_c^+(E)$ denotes the
nonnegative continuous functions with compact support,
\[
   \E\sum_{x_1,\ldots,x_m\in\Xi_n}^{\ne}
      \prod_{i=1}^mf_i(x_i)
   \longrightarrow
   \prod_{i=1}^m\int_Ef_i\,\dd\nu.
\]
Here $\sum^{\ne}$ is over ordered tuples of distinct atom indices, with
multiplicities counted in the usual factorial-measure sense.
If, for every relatively compact Borel set $B\subset E$, there is a constant $c_B>2$
such that $\sup_n\E c_B^{\Xi_n(B)}<\infty$, then
$\Xi_n$ converges in distribution, for the vague topology, to a Poisson
point process with intensity $\nu$.
\end{proposition}

This criterion follows from the factorial inversion formula
\cite[Theorem~D.6]{Betsch}, applied to independent thinning with retention
function $1-e^{-f}$, where $f\in C_c^+(E)$, and the Laplace-functional
convergence theorem \cite[Theorem~23.16]{Kallenberg}. Explicitly, with
$h=1-e^{-f}$, by factorial inversion we get that
\[
   \E e^{-\int_E f\,\dd\Xi_n}
      =\sum_{m=0}^{\infty}\frac{(-1)^m}{m!}
         \E\sum_{x_1,\ldots,x_m\in\Xi_n}^{\ne}
               \prod_{i=1}^m h(x_i),
\]
where the term with $m=0$ is one. For a relatively compact Borel set
$B\supset\operatorname{supp}f$ and integer $M\ge0$, the absolute tail after
order $M$ is at most
$\E[2^{\Xi_n(B)}\1_{\{\Xi_n(B)>M\}}]$. Uniform passage to the limit in the
inversion series is therefore justified by
\[
   \sup_n\E\!\left[2^{\Xi_n(B)}\1_{\{\Xi_n(B)>M\}}\right]
   \le\left(\frac2{c_B}\right)^{M+1}
       \sup_n\E c_B^{\Xi_n(B)}\longrightarrow0,
   \qquad M\to\infty.
\]
The limiting series is
\[
   \sum_{m=0}^{\infty}\frac{(-1)^m}{m!}
       \left(\int_E h\,\dd\nu\right)^m
   =\exp\left\{-\int_E(1-e^{-f})\,\dd\nu\right\},
\]
which is the Laplace functional of the stated Poisson process.


\renewcommand{\bibliofont}{\fontsize{10}{11}\selectfont}

\end{document}